\documentclass[12]{amsart}

\newcommand{\Z}{{\mathbb Z}}

\renewcommand{\d}{\partial}
\newcommand{\dbar}{\bar\partial}

\newcommand{\E}{{\mathbf E}}

\renewcommand{\phi}{\varphi}

\newcommand{\al}{\alpha}

\newcommand{\la}{\lambda}
\newcommand{\ep}{\varepsilon}

\newcommand{\om}{\omega}

\newcommand{\half}{{\frac{1}{2}}}

\renewcommand{\phi}{\varphi}

\newtheorem{theo}{{\sc Theorem}}[section]
\newtheorem{cor}[theo]{{\sc Corollary}}
\newtheorem{rem}[theo]{{\sc Remark}}

\newtheorem{lem}[theo]{{\sc Lemma}}

\newtheorem{prop}[theo]{{\sc Proposition}}

\title[Inverse Spectral Problems for Schr\"odinger operators ]{Inverse Spectral Problems for Schr\"odinger operators }

\author{HAMID HEZARI}

\address{Department of Mathematics, Johns Hopkins University, Baltimore, MD
21218, USA}

\email{hhezari@math.jhu.edu}

\date{November, 2007}

\begin{document}

\maketitle

\begin{abstract} In this article we improve some of
the inverse spectral results proved by Guillemin and Uribe in
\cite{GU}. They proved that under some symmetry assumptions on the
potential $V(x)$, the Taylor expansion of $V(x)$ near a
non-degenerate global minimum can be recovered from the knowledge
of the low-lying eigenvalues of the associated Schr\"odinger
operator in $\mathbb R^n$. We prove some similar inverse spectral
results using fewer symmetry assumptions. We also show that in
dimension 1, no symmetry assumption is needed to recover the
Taylor coefficients of $V(x)$. We establish our results by finding
some explicit formulas for wave invariants at the bottom of the
well.


\end{abstract}

\section{Introduction and Statement of Results}
In this article we study some inverse spectral problems of the
eigenvalue problem for the semi-classical Schr\"odinger operator,

\begin{equation} \label{Hhat} \hat P=-\frac{{\hbar}^2}{2}\triangle+V(x)
 \qquad \text{on} \quad L^2(\mathbb{R}^n),
\end{equation}
associated to the Hamiltonian
$$ P(x,\xi)=\half \,\xi^2+ V(x).$$
Here the potential $V(x)$ in ($\ref{Hhat}$) satisfies
\begin{equation} \label{Hypothesis}
 \left\{\begin{array}{lll}
 V(x) \in C^{\infty}(\mathbb{R}^n), \\
  V(x) \,\, \text{has a unique non-degenerate global minimum at} \,\, x=0 \,\, \text{and} \,\,V(0)=0, \\
  \text{For some} \,\, \varepsilon >0,
\,\,V^{-1}[0,\varepsilon]\,\, \text{is compact}. \\
\end{array} \right.
\end{equation}
Under these conditions for sufficiently small $\hbar$, say $\hbar
\in (0,h_0)$, and sufficiently small $\delta$, a classical fact
tells us the spectrum of $\hat P$ in the energy interval $[0,
\delta]$ is finite. We denote these eigenvalues by
$$ \{E_j(\hbar)\}_{j=0}^{m}.$$
We call these eigenvalues the low-lying eigenvalues of $\hat P$.
We notice the Weyl's law reads
\begin{equation} \label{hWeyl} m=N_{\hbar}(\delta)=\sharp \{j; 0 \leq E_j(\hbar)\leq
\delta\}=\frac{1}{({2\pi\hbar})^n}(\int_{\half\,\xi^2+V(x)\leq
\delta} dxd\xi + o(1)).
\end{equation}\\

Recently in \cite{GU}, Guillemin and Uribe raised the question
whether we can recover the Taylor coefficients of $V$ at $x=0$
from the low-lying eigenvalues $E_j(\hbar)$. They also established
that if we assume some symmetry conditions on $V$, namely
$V(x)=f(x_1^2,...,x_n^2)$, then the $1$-parameter family of
low-lying eigenvalues, $\{E_j(\hbar)\,|\, \hbar \in (0,h_0)\}$,
determines the Taylor coefficients of $V$ at $x=0$.\\

In this article we will attempt to recover as much of $V$ as
possible from the family $E_j(\hbar)$, by establishing some new
formulas for the wave invariants at the bottom of the potential
(Theorem $\ref{Theorem1}$). Using these new expressions for the
wave invariants, in Theorem $\ref{Theorem2}$ we improve the
inverse spectral results of \cite{GU} for a larger class of potentials.\\

A classical approach in studying this problem is to examine the
asymptotic behavior as $\hbar \to 0$ of the truncated trace
\begin{equation} \label{ttrace}
Tr(\Theta(\hat P)e^{\frac{-it}{\hbar} \hat P}),
\end{equation}
where $\Theta \in C_{0}^{\infty}([0,\infty))$ is supported in
$I=[0, \delta]$ and equals one in a neighborhood of $0$.\\

The asymptotic behavior of the truncated trace $Tr(t,\hbar)$
around the equilibrium point $(x, \xi)=(0,0)$ has been extensively
studied in the literature. It is known that (see for example
\cite{BPU}) for $t$ in a sufficiently small interval $(0,t_0)$,
$Tr(\Theta(\hat P)e^{\frac{-it}{\hbar} \hat P})$ has an asymptotic
expansion of the following form:

\begin{equation} \label{waveinvariants} Tr(\Theta(\hat P)e^{\frac{-it}{\hbar} \hat P}) \sim \sum_{j=0}^{\infty}
a_j(t)\hbar^j, \qquad \hbar \to 0.
\end{equation}
Throughout this paper when we refer to wave invariants at the
bottom of the well, we mean the
coefficients $a_j(t)$ in ($\ref{waveinvariants}$).\\

By applying an orthogonal change of variable, we can assume that
$V$ is of the form
\begin{equation}\label{W} V(x)=\half
\sum_{k=1}^{n}\omega_k^2 x_k^2+W(x), \qquad \om_k>0,\end{equation}
$$ W(x)=O(|x|^3), \qquad |x|\to 0.$$
In addition to conditions in $(\ref{Hypothesis})$, we also assume
that $\{\omega_k\}$ are linearly independent over $\mathbb Q$. We
note that we have $W(0)=\nabla W(0)=\text{Hess}W(0)=0$.
\\

Our first result finds explicit formulas for the wave invariants.
\begin{theo}\label{Theorem1} For $0< t < \min_{1\leq k \leq n}
\{\frac{\pi}{2w_k}\}$, \item[1.] $$a_0(t)=Tr(e^{\frac{-it}{\hbar}
\hat H_0})=\prod_{k=1}^n \frac{1}{2i\sin \frac{\om_k t}{2}},
\qquad \text{where} \qquad \hat H_0=-\half \hbar^2 \Delta + \half
\sum_{k=1}^{n}\omega_k^2 x_k^2.$$
\\
\item[2.] For $j\geq1$ , the wave invariants $a_j(t)$ defined in
(\ref{waveinvariants}) are given by

\begin{equation}\label{ajt}
a_j(t)=a_0(t)\sum_{l=1}^{2j}i^{l(n-1)+\frac{n}{2}}e^{{i\frac{\pi}{4}\text{sgn}H_l}}\int_0^t\int_0^{s_1}...\int_{0}^{s_{l-1}}
P_{l+j}b_l(0)ds_l...ds_1,
\end{equation}
where for every $m$,
$$P_mb_l(0)=\frac{i^{-m}}{2^m\,m!}<H_l^{-1}\nabla,\nabla>^m(b_l)(0),$$
$$ b_l=\prod_{i=1}^l W(\frac
{\cos\om_ks_i}{2}(z_{i+1}^k+z_i^k)-\frac{\sin\om_ks_i}{\om_k}\xi_i^k+(\frac{\sin\om_k(t-s_i)+\sin\om_ks_i}{\sin\om_kt})x^k),$$
and $H_l^{-1}$ is the inverse matrix of the Hessian $H_l=$Hess
$\Psi_l(0)$, where
$$\Psi_l=\Psi_l(t,x,z_1,...,z_l,\xi_1,...,\xi_l)=\sum_{k=1}^n\{(-\om_k\tan
\frac{\om_k}{2}t)x_k^2+(\frac{\om_k}{2}\cot
\om_kt)(z_1^k)^2+\sum_{i=1}^l(z_{i+1}^k-z_i^k)\xi_i^k\}.$$ The
Hessian of $\Psi_l$ is calculated with respect to every variable
except $t$. Therefore the entries of the matrix $H_l^{-1}$ are
functions of $t$. The matrix $H_l^{-1}$ is shown in
$(\ref{H_linverse})$.
\\

\item[3.] The wave invariant $a_j(t)$ is a polynomial of degree
$2j$ of the Taylor coefficients of $V$. The highest order of
derivatives appearing in $a_j(t)$ are of order $2j+2$. In fact
these higher order derivatives appear in the linear term of the
polynomial and
\begin{equation} \label{linearterm}
a_j(t)=\frac{a_0(t)}{(2i)^{j+1}}\sum_{|\vec \al|=j+1}
\frac{t}{\vec \al !}(\frac{-1}{2\vec \om}\cot\frac{\vec
\om}{2}t)^{\vec \al}D_{2 \vec \al}^{2j+2}V(0)+\{\text{a polynomial
of Taylor coefficients of order $\leq 2j+1$}\}
\end{equation}
\end{theo}
Notice that in (\ref{linearterm}), we have used the standard
shorthand notations for multi-indices, i.e. $\vec
\al=(\al_1,...\al_n)$, $\vec \om=(\om_1,...\om_n)$, $|\vec
\al|=\al_1+...+\al_n$, $\vec \al !=\al_1!...\al_n!$, $\vec X
^{\vec \al}=X_1^{\al_1}...X_n^{\al_n}$, and $D_{ \vec
\al}^{m}=\frac{\d^m }{\d x_1 ^{\al_1}...\d x_n ^{\al_n} }$ with
$m=|\vec \al|$.
\\

Our second result improves the result of Guillemin and Uribe in
\cite{GU}. This theorem is actually a non-trivial corollary of
Theorem \ref{Theorem1}.
\begin{theo} \label{Theorem2} Let $V$ satisfy $(\ref{Hypothesis})$, $(\ref{W})$,  and be of the form
\begin{equation} \label{fg} V(x)=f(x_1^2,...,x_n^2)+ x_n^3g(x_1^2,...,x_n^2),  \end{equation}
for some $f,g\in C^{\infty}(\mathbb R ^n)$. Then the low-lying
eigenvalues of $\hat P=-\half \hbar ^2 \Delta+V$ determine
$D^{|\vec \al|}_{\vec \al}V(0)$, $|\vec \al|=2,3$, and if
$D^3_{3\vec e_n}V(0):=\frac{\d^3V}{\d x_n^3}(0)\neq 0$, they
determine all the Taylor coefficients of $V$ at $x=0$.
\end{theo}
One quick consequence of Theorem $\ref{Theorem2}$ is the
following:
\begin{cor} \label{corollary} If $n=1$, and $V\in C^{\infty} (\mathbb R)$ satisfies $(\ref{Hypothesis})$,
then (with no symmetry assumptions) the low-lying eigenvalues
determine $V''(0)$ and $V^{(3)}(0)$, and if $V^{(3)}(0)\neq 0$,
then these eigenvalues determine all the Taylor coefficients of
$V$ at $x=0$.
\end{cor}
Let us briefly sketch our main ideas for the proofs. First,
because of a technical reason which arises in the proofs, we will
need to replace the Hamiltonian $P$ by the following Hamiltonian
$H$

$$ \left\{\begin{array}{lll}
  H(x,\xi)=\half \,\xi^2+V_\hbar(x), \\ \\
  V_\hbar(x)=\half \sum_{k=1}^{n}\omega_k^2
x_k^2+W_\hbar(x),\\ \\
W_\hbar(x)=\chi(\frac{x}{\hbar^{\half -\ep}})W(x),  \qquad \ep>0
\;\; \text{sufficiently small},
 \\ \end{array} \right.$$
where the cut off $\chi \in C_0^{\infty}(\mathbb R^n)$ is
supported in the unit ball $B_1(0)$ and equals one in
$B_{\half}(0).$

Then in two lemmas \,(Lemma \ref{criticallem} and Lemma
\ref{reductionlemma}) we show that for $t$ is a sufficiently small
interval $(0,t_0)$, in the sense of tempered distributions we have

$$Tr(\Theta(\hat P)e^{\frac{-it}{\hbar}\hat
P})=Tr(e^{\frac{-it}{\hbar}\hat H})+O(\hbar^\infty).$$ This
reduces the problem to studying the asymptotic of
$Tr(e^{\frac{-it}{\hbar}\hat H})$. For this we use the
construction of the kernel $k(t,x,y)$ of the propagator $U(t)=
e^{\frac{-it}{\hbar} \hat H}$ found in \cite{Z}. In \cite{Z} it is
shown that

\begin{equation} \label{k(t,x,y)introduction} k(t,x,y)=C(t)e^{\frac{i}{\hbar}S(t,x,y)}\sum_{l=0}^{\infty}
a_l(t,\hbar,x,y), \end{equation} where \begin{equation}
\label{harmonicphase} S(t,x,y)=\sum_{k=1}^{n}\frac{\om_k}{\sin
\om_k t}(\half (\cos \om_k t)(x_k^2+y_k^2)-x_ky_k),\end{equation}
$a_0=1$, and for $l \geq 1$,

$$a_l(t,\hbar,x,y)=(\frac{-1}{2\pi})^{ln}(\frac{1}{i\hbar})^{l(n+1)}
\int_0^t...\int_0^{s_{l-1}} \overbrace{\int ...
\int}^{2l}e^{\frac{i}{\hbar}\Phi_l} b_l(s,x,y,\vec z,\vec
\xi)d^lzd^l \xi d^ls,$$ where
$$\Phi_l=\sum_{k=1}^n\{\frac{\om_k}{2}\cot
\om_kt(z_1^k)^2+\sum_{i=1}^l(z_{i+1}^k-z_i^k)\xi_i^k\},
$$ and $$b_l=\prod_{i=1}^l W_\hbar(\frac
{\cos\om_ks_i}{2}(z_{i+1}^k+z_i^k)-\frac{\sin\om_ks_i}{\om_k}\xi_i^k+\frac{\sin\om_k(t-s_i)}{\sin\om_kt}y^k+
\frac{\sin\om_ks_i}{\sin\om_kt}x^k).
$$\\

Next we apply the expression in (\ref{k(t,x,y)introduction}) for
$k(t,x,y)$ to the formula $Tr(e^{\frac{-it}{\hbar}\hat H})=\int
k(t,x,x)dx$. Then we obtain an infinite series of oscillatory
integrals, each one corresponding to one $a_l$. Finally we apply
the method of stationary phase to each oscillatory integral and we
show that the resulting series is a valid asymptotic expansion.
From the resulting asymptotic expansion we obtain the formulas
($\ref{ajt}$).
\\

Now let us compare our approach for the construction of $k(t,x,y)$
with the classical approach. In the classical approach (see for
instance \cite{DSj},\cite{D}, \cite{R}, \cite{BPU} and \cite{U}),
one constructs a WKB approximation for the kernel $k_P(t,x,y)$ of
the operator $\Theta(\hat P)e^{\frac{-it}{\hbar} \hat P}$, i.e.
\begin{equation}\label{WKB}
k_P(t,x,y)=\int e^{\frac{i}{\hbar}(\phi_P(t,x,\eta)-y.\eta)}
b_P(t,x,y,\eta,\hbar) d\eta,
\end{equation}
where $\phi_P(t,x,\eta)$ satisfies the Hamilton-Jacobi equation
(or eikonal equation in geometrical optics)

$$ \d_t\phi_P(t,x,\eta)+P(x,\d_x\phi_P(t,x,\eta))=0, \qquad
\phi_P|_{t=0}=x.\eta, $$
and the function $b_P$ has an asymptotic
expansion of the form
$$b_P(t,x,y,\eta,\hbar) \sim \sum_{j=0}^{\infty} b_{P,j}(t,x,y,\eta) \hbar ^j. $$
The functions $b_{P,j}(t,x,y,\eta)$ are calculated from the so
called transport equations. See for example \cite{R}, \cite{DSj},
\cite{EZ} or Appendix A of the paper in hand for the details of
the above construction.

In this setting, when one integrates the kernel $k(t,x,y)$ on the
diagonal and applies the stationary phase to the given oscillatory
integral, one obtains very complicated expressions for the wave
invariants. Of course the classical calculations above show the
existence of asymptotic formulas of the form
(\ref{waveinvariants}) (which can be used to get Weyl-type
estimates for the counting functions of the eigenvalues, see for
example \cite{BPU}). Unfortunately these formulas for the wave
invariants are not helpful when trying to establish some inverse
spectral
results.\\

Hence, one should look for more efficient methods to calculate the
wave invariants $a_j(t)$. One approach is to use the
semi-classical Birkhoff normal forms, which was used in \cite{GU},
to get some inverse spectral results as mentioned in the beginning
of the introduction. The Birkhoff normal forms methods were also
used before by S. Zelditch in \cite{Z4} to obtain positive inverse
spectral results for real analytic domains with symmetries of an
ellipse. Zelditch proved that for a real analytic plane domain
with symmetries of an ellipse, the wave invariants at a bouncing
ball orbit, which is preserved by the symmetries, determine the
real analytic domain under isometries of the domain.
\\

Recently in \cite{Z3}, Zelditch improved his earlier result to
real analytic domains with only one mirror symmetry. His approach
for this new result was different. He used a direct approach
(Balian-Bloch trace formula) which involves Feynman-diagrammatic
calculations of the stationary phase method to obtain a more
explicit formula for the wave invariants at the bouncing ball
orbit.

Motivated by the work of Zelditch \cite{Z3} mentioned above, our
approach in this article is also somehow direct and involves
combinatorial calculations of the stationary phase.

Our formula in $(\ref{k(t,x,y)introduction})$ for the kernel of
the propagator, $U(t)= e^{\frac{-it}{\hbar} \hat H}$, is different
from the WKB-expression in the sense that we only keep the
quadratic part of the phase function, namely the phase function
$S(t,x,y)$ in $(\ref{harmonicphase})$ of the propagator of
Anisotropic oscillator , and we put the rest in the amplitude
$\sum_{l=0}^{\infty} a_l(t,\hbar,x,y)$ in
$(\ref{k(t,x,y)introduction})$. The details of this construction
are mentioned in Section $2.2$.\\

\begin{rem} After the initial posting of this article,  Guillemin and Colin de Verdi\`{e}re
posted two articles (see \cite{CG1}, also \cite{C}) in which they
study some inverse spectral problems of 1 dimensional
semi-classical Schr\"odinger operators. One of the main results in
\cite{CG1} is our Corollary \ref{corollary} in this paper.
\end{rem}

\subsection{Acknowledgements:} I am sincerely grateful to Steve
Zelditch for introducing the problem and many helpful discussions
and suggestions on the subject. I would also like to thank him for
his great support and encouragement as I was writing this article.

\section{Proofs of the results}

\subsection{Some reductions}

Because of some technical issues arising in the proof of Theorem
\ref{Theorem1}, we will need to use the following lemmas as
reductions.
\\

In the following, we let $\chi \in C_0^{\infty}(\mathbb R^n)$ be a
cut off which is supported in the unit ball $B_1(0)$ and equals
one in $B_{\half}(0)$.

\begin{lem}\label{criticallem} Let the Hamiltonians $P$ and $H$ be
defined by \begin{equation} \label{PH} \left\{ \begin{array}{ll}
  P(x,\xi)=\half \,\xi^2+ V(x) \\ \\
  V(x)=\half\sum_{k=1}^{n}\omega_k^2 x_k^2+W(x) \\
\end{array}\right., \qquad
\left\{\begin{array}{lll}
  H(x,\xi)=\half \,\xi^2+V_\hbar(x), \\ \\
  V_\hbar(x)=\half \sum_{k=1}^{n}\omega_k^2
x_k^2+W_\hbar(x),\\ \\
W_\hbar(x)=\chi(\frac{x}{\hbar^{\half -\ep}})W(x),  \qquad \ep>0
\;\; \text{sufficiently small},
 \\ \end{array} \right.
\end{equation}
and let $\hat P$ and $\hat H$ be the corresponding Weyl (or
standard) quantizations. Then for $t$ in a sufficiently small
interval $(0,t_0)$

$$ Tr(\Theta(\hat P)e^{\frac{-it}{\hbar} \hat P})=Tr(\Theta(\hat H)e^{\frac{-it}{\hbar} \hat
H})+O(\hbar ^\infty).$$ In other words, the wave invariants
$a_j(t)$ will not change if we replace $P$ by $H$.
\end{lem}

\textsc{Proof.} Proof is given in Appendix A.
\\

Next we use the following lemma to get rid of $\Theta(\hat H)$.

\begin{lem} \label{reductionlemma} Let $H$ be defined by
$(\ref{PH})$. Then in the sense of tempered distributions
$$ Tr(\Theta(\hat H)e^{\frac{-it}{\hbar} \hat H})=Tr(e^{\frac{-it}{\hbar} \hat
H})+O(\hbar ^\infty).$$ This means that if we sort the spectrum of
$\hat H$ as

$$ E_1(\hbar)<E_2(\hbar) \leq...\leq E_j(\hbar) \to +\infty,$$
then for every Schwartz function $\phi(t) \in S(\mathbb R)$ $$
<Tr(e^{\frac{-it}{\hbar} \hat H})-Tr(\Theta(\hat
H)e^{\frac{-it}{\hbar} \hat H}), \,\phi(t)>=\sum_{j=1}^{\infty}
(1-\Theta(E_j(\hbar)))\hat
\phi(\frac{E_j(\hbar)}{\hbar})=O(\hbar^{\infty}).$$
\end{lem}

 \textsc{Proof.} Proof is given in Appendix B.
 \\

Because of the above lemmas, it is enough to study the asymptotic
of $Tr(e^{\frac{-it}{\hbar} \hat H})$.

\subsection{Construction of $k(t,x,y)$, the kernel of $e^{\frac{-it}{\hbar} \hat
H}$}

In this section we follow the construction in \cite{Z} to obtain
an oscillatory integral representation of $k(t,x,y)$, the kernel
of the propagator $e^{\frac{-it}{\hbar} \hat H}$. The reader
should consult \cite{Z} for many details. In that article Zelditch
uses the Dyson's Expansion of propagator to study the
singularities of the kernel $k(t,x,y)$. But he does not consider
the semi-classical setting $\hbar\to 0$ in his calculations
\,(i.e. $\hbar=1$). So we follow the same calculations but also
consider $\hbar$ carefully.

In \cite{Z}, potentials of the form

$$V(x)=\half
\sum_{k=1}^{n}\omega_k^2 x_k^2+ W_\hbar(x),$$ are considered,
where $W_\hbar \in B(\mathbb R^n)$, i.e. bounded with bounded
derivatives.
\\

We denote
$$\left\{\begin{array}{ll}
  \hat H_0=-\half \hbar^2 \Delta + \half \sum_{k=1}^{n}\omega_k^2
x_k^2, \qquad \text{(Anisotropic Oscillator)} \\ \\
  \hat H=\hat H_0+W_\hbar(x)=-\half \hbar^2 \Delta + V_\hbar(x), \\
\end{array}\right.$$
and by $U_0(t)=e^{\frac{-it}{\hbar} \hat H_0}$, and
$U(t)=e^{\frac{-it}{\hbar} \hat H}$, we mean their corresponding
propagators.
\\

From $$(i\hbar \d_t-\hat H_0)U(t)=W_\hbar.U(t),$$ we obtain
\begin{equation} \label{Duhamel}U(t)=U_0(t)+\frac{1}{i\hbar}\int_0^t U_0(t-s).W_\hbar. U(s)ds. \end{equation}
By iteration we get the norm convergent Dyson Expansion:

\begin{equation} \label{Dyson}U(t)=U_0(t)+\sum_{l=1}^{\infty}
\frac{1}{(i\hbar)^l}\int_0^t...\int_0^{s_{l-1}} U_0(t)
[U_0(s_1)^{-1}.W_\hbar. U_0(s_1)]...[U_0(s_l)^{-1}. W_\hbar.
U_0(s_l)]ds_{l}...ds_1.
\end{equation}
For $t \neq \frac{m\pi}{\om_k}$, the kernel of $U_0(t)$ is given
by

\begin{equation} \label{k_0} k_0(t,x,y)=(\prod_{k=1}^n
\frac{\om_k}{2\pi i \hbar \sin \om_k
t})^{\half}e^{\frac{i}{\hbar}S(t,x,y)}, \end{equation} where
$$S(t,x,y)=\sum_{k=1}^{n}\frac{\om_k}{\sin \om_k t}(\half (\cos
\om_k t)(x_k^2+y_k^2)-x_ky_k) .$$ Then by taking kernels in
$(\ref{Dyson})$ and after some change of variables (see \cite{Z}),
we get

\begin{equation} \label{k(t,x,y)} k(t,x,y)=(\prod_{k=1}^n
\frac{\om_k}{2\pi i \hbar \sin \om_k
t})^{\half}e^{\frac{i}{\hbar}S(t,x,y)}\sum_{l=0}^{\infty}
a_l(t,\hbar,x,y),
\end{equation}
where $a_0=1$ and for $l \geq 1$,
\begin{equation} \label{a_l} a_l(t,\hbar,x,y)=(\frac{-1}{2\pi})^{ln}(\frac{1}{i\hbar})^{l(n+1)}
\int_0^t...\int_0^{s_{l-1}} \overbrace{\int ...
\int}^{2l}e^{\frac{i}{\hbar}\Phi_l} b_l(s,x,y,\vec z,\vec
\xi)d^lzd^l \xi d^ls,
\end{equation}
where \begin{equation} \label{Phi_l}
\Phi_l=\sum_{k=1}^n\{\frac{\om_k}{2}\cot
\om_kt(z_1^k)^2+\sum_{i=1}^l(z_{i+1}^k-z_i^k)\xi_i^k\},
\end{equation}
and \begin{equation} \label{b_l} b_l=\prod_{i=1}^l W_\hbar(\frac
{\cos\om_ks_i}{2}(z_{i+1}^k+z_i^k)-\frac{\sin\om_ks_i}{\om_k}\xi_i^k+\frac{\sin\om_k(t-s_i)}{\sin\om_kt}y^k+
\frac{\sin\om_ks_i}{\sin\om_kt}x^k).\qquad \quad(z_{l+1}:=0)
\end{equation}

In \cite{Z}, integration by parts in the variables $z,\xi,$ for
the integrals in $(\ref{a_l})$ are performed (in that article,
there is no $\hbar$, i.e. $\hbar=1$, because of different
motivations)\,to prove that the sum
$a(t,1,x,y)=\sum_{l=0}a_l(t,1,x,y)$ is absolutely uniformly
convergent. Additionally, it is shown that $a(t,1,x,y)$ is in
$B(\mathbb R_x ^n \times \mathbb R_x ^n).$ More precisely, there
exists $k_0=k_0(\al,\beta,n)$ such that

\begin{equation} \label{daldbeta}
|\d^{\vec \al}_x \d^{\vec \beta}_y a_l(t,1,x,y)| \leq \frac{1}{l!}
 {C_{\al,\beta,n}(t)}^l||W_1||_{|\al|+|\beta|+k_0}^l, \qquad \quad
 (W_1=W_\hbar\mid_{\hbar=1})
\end{equation}
which implies
$$|\d^{\vec \al}_x \d^{\vec \beta}_y a(t,1,x,y)| \leq
\exp\{{C_{\al,\beta,n}(t)||W_1||_{|\al|+|\beta|+k_0}}\}.$$ The
estimates ($\ref{daldbeta}$) will change if one considers $\hbar$
in the calculations. We will establish these $\hbar$ estimates in
Lemma $\ref{lem2}$. As a simple consequence of Lemma $\ref{lem2}$,
let us assume for now that $a(t,\hbar,x,y) \in B(\mathbb R_x ^n
\times \mathbb R_x ^n).$

\subsection{Oscillatory Integral Representation for the Trace of $U(t)=e^{-\frac{it}{\hbar} \hat H}$}
In this section we show that the integral $Tr\,U(t)=\int
k(t,x,x)dx$ is convergent as an oscillatory integral and that from
($\ref{k(t,x,y)}$) and ($\ref{a_l}$) we can write

\begin{equation} \label{Trace}
Tr\,U(t)=(\prod_{k=1}^n \frac{\om_k}{2\pi i \hbar \sin \om_k
t})^{\half}\sum_{l=0}^{\infty}(\frac{-1}{2\pi})^{ln}(\frac{1}{i\hbar})^{l(n+1)}
\int_0^t...\int_0^{s_{l-1}} \overbrace{\int ...
\int}^{2l+1}e^{\frac{i}{\hbar}\Psi_l} b_l(s,x,x,\vec z,\vec
\xi)d^ls\,d^lz\,d^l\xi\,dx,
\end{equation}
where \begin{equation} \label{Psi} \Psi_l=S(t,x,x)+\Phi_l
=\sum_{k=1}^n\{(-\om_k\tan\frac{\om_k}{2}t)\,x_k^2+
\frac{\om_k}{2}\cot
\om_kt(z_1^k)^2+\sum_{i=1}^l(z_{i+1}^k-z_i^k)\xi_i^k\}.
\end{equation}
Before proving $(\ref{Trace})$, we review some standard facts.
First of all we know that the sum
$$ Tr\,U(t)=\sum e^{-\frac{it\E_j(\hbar)}{\hbar}}$$
is convergent in the sense of tempered distributions, i.e.
$Tr\,U(t) \in S'(\mathbb R)$. This can be shown by the Weyl's law
in its high energy setting, which implies that for potentials of
the form $V(x)=\half \sum_{k=1}^{n}\omega_k^2 x_k^2+W_\hbar(x)$,
with $W\in B(\mathbb R^n)$, for fixed $\hbar$, the $j-$th
eigenvalue $E_j(\hbar)$ satisfies

\begin{equation} \label{Weyl'slaw}E_j(\hbar) \sim C(n,\hbar)j^{\frac{2}{n}}, \qquad j\to
\infty. \end{equation}
Another way to define $Tr\,U(t)$ is to
write it as the limit

\begin{equation}\label{Traceepsilon}
Tr\,U(t)=\lim_{\ep \to 0^+}Tr\,U(t-i\ep)=\lim_{\ep \to 0^+}\sum
e^{-\frac{(it+\ep)\E_j(\hbar)}{\hbar}}.
\end{equation}
This time the Weyl's law $(\ref{Weyl'slaw})$ implies that the sum
$Tr\,U(t-i\ep)$ is absolutely uniformly convergent because of the
rapidly decaying factor $e^{-\frac{\ep\E_j(\hbar)}{\hbar}}$. As a
result, $U(t-i\ep)$ is a trace class operator. It is clear that
the kernel of $U(t-i\ep)$ is $k(t-i\ep,x,y)$, the analytic
continuation of the kernel $k(t,x,y)$ of $U(t)$. Clearly
$k(t-i\ep,x,y)$ is continuous in $x$ and $y$. So we can write
$Tr\,U(t-i\ep)= \int k(t-i\ep,x,x)dx$. We notice that this
integral is uniformly convergent. This is because up to a constant
this integral equals to $\int
e^{\frac{i}{\hbar}S(t-i\ep,x,x)}\,a(t-i\ep,\hbar,x,x),$ and the
exponential factor in the integral is rapidly decaying for $\ep
>0$ as $|x| \to \infty$ and $a$ is a bounded function. More precisely
$$\Re(iS(t-i\ep,x,x))=\sum_{k=1}^n \Re(-i {\om_k}\tan(\frac{\om_k (t-i\ep)}{2}))x_k^2
=\sum_{k=1}^n\frac{\om_k(1-e^{2\ep
\om_k})}{|1+e^{w_k(it+\ep)}|^2}x_k^2,
$$
and $$\frac{\om_k(1-e^{2\ep \om_k})}{|1+e^{w_k(it+\ep)}|^2}<0.$$
The discussion above shows that the integral $\int k(t,x,x)dx$ can
be defined by integrations by parts as follows:
Since

$$ <D_x>^2
e^{iS(t,x,x)}:=(1-\Delta)e^{iS(t,x,x)}=(1+\parallel2\vec\om\tan(\frac{\vec\om
t}{2})\vec x\parallel ^2+2i\sum_{k=1}^n \om_k\tan(\frac{\om_k
t}{2}))e^{iS(t,x,x)},$$ we can write
\begin{equation} \label{intbyparts} \int
e^{\frac{i}{\hbar}S(t,x,x)}\,a(t,\hbar,x,x)dx=\hbar^{\frac{n}{2}}\int
e^{iS(t,x,x)}\,
a(t,\hbar,\sqrt{\hbar}x,\sqrt{\hbar}x)dx=\end{equation}
$$\hbar^{\frac{n}{2}}\int e^{iS(t,x,x)}\,
(<D_x>^2\,(1+\parallel2\vec\om\tan(\frac{\vec\om t}{2})\vec
x\parallel ^2+2i\sum_{k=1}^n \om_k\tan(\frac{\om_k
t}{2}))^{-1})^{n_0} a(t,\hbar,\sqrt{\hbar}x,\sqrt{\hbar}x)dx$$ If
we assume $0< t < \min_{1\leq k \leq n} \{\frac{\pi}{2w_k}\}$,
then by choosing $n_0>\frac{n}{2}$, and because $a(t,\hbar,x,y)
\in B(\mathbb R_x ^n \times \mathbb R_y ^n)$,  the integral
becomes absolutely convergent. Finally, since by
(\ref{boundfora_l}) the series
$a(t,\hbar,x,y)=\sum_{l=0}^{\infty}a_l(t,\hbar,x,y)$ is absolutely
uniformly convergent, we have $$\int
e^{\frac{i}{\hbar}S(t,x,x)}\,a(t,\hbar,x,x)dx=\sum_{l=0}^{\infty}\int
e^{\frac{i}{\hbar}S(t,x,x)}\,a_l(t,\hbar,x,x)dx,$$ and therefore
we obtain $(\ref{Trace})$, which is an infinite sum of oscillatory
integrals. The next step is to apply the stationary phase method
to each integral in $(\ref{Trace})$ and then add the asymptotics
to obtain an asymptotic expansion for the $Tr\,U(t)$. Because we
have an infinite sum of asymptotic expansions, we have to
establish that the resulting asymptotic for the trace is a valid
approximation. Hence we have to find some appropriate
$\hbar$-estimates for the remainder term of the series
$(\ref{Trace})$.

\subsection{$\hbar$-estimates for the remainder}
The goal of this section is to find an $\hbar$-estimate for the
remainder term of the series $(\ref{Trace})$. First we state the
following important estimate, which shows how the estimates in
$(\ref{daldbeta})$ change in the $\hbar$-dependence case.

\begin{lem} \label{lem2}
For every $\vec \al$, $\vec \beta$, there exists $k_0=k_0(\al,
\beta)$ such that for every $0<\hbar\leq h_0 \leq 1$

\begin{equation} \label{boundfora_l}
|\d^{\vec \al}_x \d^{\vec \beta}_y a_l(t,\hbar,x,y)| \leq
\frac{{C_{\al,\beta,n}(t)}^l||W_\hbar^\star||_{|\al|+|\beta|+k_0}^l}{l!}\hbar^{l(\half-3\ep)},
\end{equation} where
\begin{equation} \label{wstar} W_\hbar^{\star}(x)=\frac{ W_\hbar(\hbar^{\half}
x)}{\hbar^{3(\half-\ep)}}= \chi(\hbar^\ep x)\frac{W(\hbar^{\half}
x)}{\hbar^{3(\half-\ep)}}
\end{equation} is uniformly in $B(\mathbb R_x^n)$; i.e. $W_\hbar^{\star}$ is bounded with
bounded derivatives and the bounds are independent of $\hbar$.
\end{lem}
Before proving this lemma, let us show how the lemma is applied to
get estimates for the remainder term of the series
$(\ref{Trace})$. Define $I_l(t,\hbar)$ to be the $l$-th term in
$(\ref{Trace})$ (removing the constant $\prod_{k=1}^n
(\frac{\om_k}{2\pi i \sin \om_k t})^{\half}$), i.e.
\begin{equation} \label{I_l1} I_l(t,\hbar)= \hbar^{-\frac{n}{2}} \int
e^{\frac{i}{\hbar}S(t,x,x)}\,a_l(t,\hbar,x,x)dx.
\end{equation}
Hence by this notation, $Tr\,U(t)=\prod_{k=1}^n (\frac{\om_k}{2\pi
i \sin \om_k t})^{\half}\sum_{l=0}^{\infty}I_l(t,\hbar)$. If in
$(\ref{I_l1})$ we integrate by parts as we did in
$(\ref{intbyparts})$, and choose $n_0=[\frac{n}{2}]+1$, then using
$(\ref{boundfora_l})$ we get
$$ |I_l(t,\hbar)|\leq C'_n(t)
\frac{{C_{n}(t)}^l||W_\hbar^\star||_{2n_0+k_0}^l}{l!}\hbar^{l(\half-3\ep)},
$$
where $C_{n}(t)=\max_{|\al|+|\beta|\leq
2n_0}\{C_{\al,\beta,n}(t)\}.$ We choose $\ep>0$ such that
$\half-3\ep>0$, or $\ep<\frac{1}{6}$. Now it is clear that for
every positive integer $m$, and every $0<\hbar\leq h_0 \leq 1$,

\begin{equation} \label{luckyestimate}
|\sum_{l=m}^{\infty} I_l(t,\hbar)| \leq
C'_n(t)e^{\{C_{n}(t)||W_\hbar^\star||_{2n_0+k_0}\}}\hbar^{m(\half-3\ep)}.
\end{equation}
Since by Lemma $\ref{lem2}$, $\sup _{0<\hbar\leq
1}||W_\hbar^\star||_{2n_0+k_0} < \infty$, we have

\begin{equation} \label{luckyestimate2}
Tr\,U(t)= \prod_{k=1}^n (\frac{\om_k}{2\pi i \sin \om_k
t})^{\half}\sum_{l=0}^{m-1}I_l(t,\hbar) +
O(\hbar^{m(\half-3\ep)}).
\end{equation}
This is a useful estimate that we will use in the next section to
obtain the asymptotic expansion of the trace and prove Theorem
$\ref{Theorem1}$. Let us first prove Lemma $\ref{lem2}$.
\\

\noindent {\textsc{Proof of Lemma $\ref{lem2}$}:
\\
The proof is straightforward from ($\ref{daldbeta}$). First, in
$(\ref{a_l})$, we apply the change of variables $\vec z \mapsto
\hbar^{\half} \vec z$ and $\vec \xi \mapsto \hbar^{\half} \vec
\xi$. This gives us $\hbar^{ln}$ in front of the integral. Then we
replace $W_\hbar$ by $\hbar^{3(\half-\ep)}W_\hbar^\star$. After
collecting all the powers of $\hbar$ in front of the integral we
obtain

$$ a_l(t,\hbar,x,y)=(\frac{-1}{2\pi})^{ln}\hbar^{l(\half-3\ep)}\int_0^t...\int_0^{s_{l-1}} \overbrace{\int ...
\int}^{2l}e^{i\Phi_l}b_l^\star(s,x,y,\vec z, \vec \xi)  d^lzd^l
\xi d^ls,$$ where $$b_l^\star(s,x,y,\vec z, \vec
\xi)=\prod_{i=1}^l W_\hbar^\star(\frac
{\cos\om_ks_i}{2}(z_{i+1}^k+z_i^k)-\frac{\sin\om_ks_i}{\om_k}\xi_i^k+\frac{\sin\om_k(t-s_i)}{\sin\om_kt}y^k+
\frac{\sin\om_ks_i}{\sin\om_kt}x^k).$$ Next we apply
($\ref{daldbeta}$) to the above integral with $W_1$ replaced by
$W_\hbar^{\star}$, and we get $(\ref{boundfora_l})$. To finish the
proof we have to show that for every positive integer $m$ we can
find uniform bounds (i. e. independent of $\hbar$) for the $m$-th
derivatives of the function $W_\hbar^\star(x)$. Since $\chi(x)$ is
supported in the unit ball, from the definition $(\ref{wstar})$ we
see that $W_\hbar^\star$ is supported in $|x|<h^{-\ep}$. So from
(\ref{wstar}) it is enough to find uniform bounds in $\hbar$ for
the $m$-th derivatives of the function $\frac{W(\hbar^{\half}
x)}{\hbar^{3(\half-\ep)}}$ in the ball $|x|<\hbar ^{-\ep}$. This
is very clear for $m \geq 3$. For $m<3$ , we use the order of
vanishing of $W(x)$ at $x=0$. Since $W(0)=\nabla
W(0)=$Hess$W(0)=0$, the order of vanishing of $W$ at $x=0$ is 3.
Therefore in the ball $|x|<\hbar ^{-\ep}$, the functions

$$\frac{W(\hbar^{\half} x)}{(\hbar^{\half}x)^3}, \frac{\d^ \al W(\hbar^{\half}
x)}{(\hbar^{\half}x)^2}, \frac{\d^\al \d^\beta W(\hbar^{\half}
x)}{\hbar^{\half}x},$$
are bounded functions with uniform bounds
in $\hbar$, and the statement follows easily for $m<3$, noting
that $|x|<\hbar^{-\ep}$.

\subsection{Stationary phase calculations}
In this section we will apply the stationary phase method to each
$I_l(t,\hbar)$ in $(\ref{luckyestimate2})$. We know
\begin{equation} \label{I_l2}
I_l(x,\hbar)=(\frac{-1}{2\pi})^{ln}(\frac{1}{i\hbar})^{l(n+1)+\frac{n}{2}}
\int_0^t...\int_0^{s_{l-1}} \overbrace{\int ...
\int}^{2l+1}e^{\frac{i}{\hbar}\Psi_l} b_l(s,x,x,\vec z,\vec
\xi)d^ls\,d^lz\,d^l\xi\,dx.
\end{equation}

It is easy to see that the only critical point of the phase
function $\Psi_l$, given by (\ref{Psi}), is at $(x,\vec z, \vec
\xi\,)=0$.
\\

Next we calculate $H_l=$Hess$\Psi_l(0)$ and $H_l^{-1}$. In the
following we use the notation $D(\vec v)$ for the diagonal matrix
$Diag(v_1,...,v_n)$, where $\vec v=(v_1,...,v_n)$. From
$(\ref{Psi})$, we get
\begin{equation} \label{H_l} H_l=\left(%
\begin{array}{ccc}
  D(-2\vec \om \tan(\frac{\vec \om t}{2}))_{n \times n} & 0 & 0 \\
  0 & \left(%
\begin{array}{cc}
  D(\vec \om \cot\vec \om t)_{n\times n} & 0 \\
  0 &  0\\
\end{array}%
\right)_{ln \times ln} & \left(%
\begin{array}{ccccc}
  \!\!-I &  0 & .. &  & 0 \\
  I & \!\!\!-I & 0 & .. & 0 \\
  0 & I & ._{._.} &  &  \\
  \dot . &  & ._{._.} &  &  \\
  0 &  &  & I & \!\!-I \\
\end{array}%
\right)_{ln \times ln} \\
  0 & \left(%
\begin{array}{ccccc}
  \!\!-I &  I & .. &  & 0 \\
  0 & \!\!\!-I & I & .. & 0 \\
  . & 0& ._{._.} & ._{._.} &  \\
  . &  &    &  & I  \\
  0 &  &  & 0 & \!\!-I \\
\end{array}%
\right)_{ln \times ln} & 0 \\
\end{array}%
\right)_{(2l+1)n\times (2l+1)n},
\end{equation}
where $I=I_{n\times n}$ is the identity matrix of size $n\times
n$.

Since $H_l$ is of the form $H_l=\left(%
\begin{array}{ccc}
  K & 0 & 0 \\
  0 & A & B \\
  0 & B^T & 0 \\
\end{array}%
\right)$, the inverse matrix equals $H_l^{-1}= \left(%
\begin{array}{ccc}
  K^{-1} & 0 & 0 \\
  0 & 0 & {B^T}^{-1} \\
  0 & B^{-1} & -B^{-1}A{B^T}^{-1} \\
\end{array}%
\right).$ A simple calculation shows that

\begin{equation} \label{H_linverse}
H_l^{-1}=\left(%
\begin{array}{ccc}
  D(\frac{-1}{2\vec \om}\cot(\frac{\vec \om t}{2})) & 0 & 0 \\
  0 & 0 & \left(%
\begin{array}{ccccc}
  -I &  -I & ... &  & -I \\
   0 & -I & -I & ... & -I \\
  0 & 0 & -I & ._{._.} &  \dot .\\
  \dot . & ._{._.} & ._{._.} &  & -I \\
  0 & ... &  & 0& -I \\
\end{array}%
\right) \\
  0 & \left(%
\begin{array}{ccccc}
  -I &  0 & ... &  & 0 \\
   -I & -I & 0 &  & \dot .\\
   \dot .& \dot . &._{._.} &  & \dot .\\
   -I& ._{._.} & ._{._.} & -I &0 \\
  -I & ... &  & -I& -I \\
\end{array}%
\right)& \left(%
\begin{array}{cccc}
  -\Omega & -\Omega & ... & -\Omega \\
  -\Omega & -\Omega & ... & -\Omega \\
\dot{\dot .} &  &  & \dot{\dot .} \\
  -\Omega & -\Omega & ... & -\Omega \\
\end{array}%
\right) \\
\end{array}%
\right),
\end{equation}
where $\Omega=D(\vec w \cot \vec \om t)$.
\\

It is also easy to see that

\begin{equation} \label{det}
\text{det}H_l=(-1)^{(l+1)n} \prod_{k=1}^n 2\om_k\tan \frac{\om_k
t}{2}.
\end{equation}
By applying the stationary phase lemma to (\ref{I_l2}) and
plugging into $(\ref{luckyestimate2})$ we obtain
\begin{equation} \label{Trace2} Tr\,U(t)=\prod_{k=1}^n (\frac{\om_k}{\sin \om_k
t})^{\half}\sum_{l=0}^{m-1}\frac{\hbar^{-l}}{{i}^{l(n+1)+\frac{n}{2}}}
\frac{(-1)^{ln}e^{i\frac{\pi}{4}\text{sgn}H_l}}{\sqrt{|\text{det}H_l|}}\int_0^t...
\int_0^{s_{l-1}}\sum_{j=0}^{\infty}\hbar^jP_jb_l(0)ds_l...ds_1+O(\hbar^{m(\half-3\ep)}),\end{equation}
where
\begin{equation} \label{P_jb_l}
P_jb_l(x,\vec z,\vec
\xi)=\frac{i^{-j}}{2^j.j!}<H_l^{-1}\nabla,\nabla>^j\,b_l(x,\vec
z,\vec \xi)=\frac{i^{-j}}{2^j.j!} \sum_{r_1,...,r_{2j}\in\mathcal
A_l}h_l^{r_1r_2}...h_l^{r_{2j-1}r_{2j}}\frac{\d ^{2j}b_l(x,\vec
z,\vec \xi)}{\d{r_1}...\d{r_{2j}}},\end{equation} where in the sum
(\ref{P_jb_l}) the indices $r_1,...,r_{2j}$ run in the set
$\mathcal A_l=\{x^k,z_1^k,..z_l^k,\xi_1^k,..\xi_l^k\}_{k=1}^n$,
and $h_l^{rr'}$ with $r,r'\in \mathcal A_l$, corresponds to the
$(r,r')$-th entry of the inverse Hessian $H_l^{-1}$.
\\

We note that $P_jb_l(0)=0$ if $2j<3l$. This is true because of
$(\ref{b_l})$ and because $W(0)=\nabla W(0)=\text{Hess}W(0)=0$.
This implies, first, there are not any negative powers of $\hbar$
in the expansion (as we were expecting). Second, the constant term
\,(i.e. the $0$-th wave invariant), which corresponds to the term
$l=j=0$ in the sum, equals

$$a_0(t)=Tr U_0(t)=\prod_{k=1}^n (\frac{\om_k}{\sin \om_k
t})^{\half}\frac{i^{-\frac{n}{2}}e^{\frac{-\pi n}{4}}}{(2\om_k\tan
\frac{\om_kt}{2})^{\half}}=\prod_{k=1}^n \frac{1}{2i\sin
\frac{\om_k t}{2}}.$$
And third (using $(\ref{det})$), for $j \geq
1$ the coefficient of $\hbar^j$ in $(\ref{Trace2})$ equals

$$a_j(t)=(\prod_{k=1}^n \frac{1}{2i\sin
\frac{\om_k
t}{2}})\sum_{l=1}^{2j}i^{l(n-1)+\frac{n}{2}}e^{{i\frac{\pi}{4}\text{sgn}H_l}}\int_0^t\int_0^{s_1}...\int_{0}^{s_{l-1}}
P_{l+j}b_l(0)ds_l...ds_1.$$ The sum goes only up to $2j$ because
if $l>2j$ then $2(l+j)<3l$ and $P_{l+j}b_l(0)=0$.
\\ \\
This proves the first two parts of Theorem $\ref{Theorem1}$.

\subsection{Calculations of the wave invariants and the proof of Theorem $\ref{Theorem1}$.3.}
In this section we try to calculate the wave invariants $a_j(t)$
from the formulas $(\ref{ajt})$. First of all, let us investigate
how the terms with highest order of derivatives appear in
$a_j(t)$. Because $b_l$ is the product of $l$ copies of $W_\hbar$
functions, and because we have to put at least $3$ derivatives on
each $W_\hbar$ to obtain non-zero terms, the highest possible
order of derivatives that can appear in $P_{j+l}b_l(0)$, is
$2(j+l)-3(l-1)=2j-l+3$. This implies that, because in the sum
$(\ref{ajt})$ we have $1\leq l \leq 2j$, the highest order of
derivatives in $a_j(t)$ is $2j+2$ and those derivatives are
produced by the term corresponding to $l=1$, i.e. $P_{j+1}b_1(0)$.
The formula $(\ref{ajt})$ also shows that $a_j(t)$ is a polynomial
of degree $2j$. The term with the highest polynomial order is the
one with $l=2j$, i.e. $P_{3j}b_{2j}(0)$ (which has the lowest
order of derivatives) and the term $P_{j+1}b_1(0)$ is the linear
term of the polynomial. Now let us calculate $P_{j+1}b_1(x,\vec
z,\vec \xi\,\,)$ and prove Theorem $\ref{Theorem1}$.3.
\\

By $(\ref{P_jb_l})$,

\begin{equation} \label{P_j+1b_1}
P_{j+1}b_1= \frac{i^{-(j+1)}}{2^{j+1}.(j+1)!}
\sum_{r_1,...,r_{2j+2}\in \mathcal
A_1}h_1^{r_1r_2}...h_1^{r_{2j+1}r_{2j+2}}\frac{\d
^{2j+2}b_1}{\d{r_1}...\d{r_{2j+2}}},
\end{equation}
where here by (\ref{b_l})

$$b_1=W_\hbar(\frac
{\cos\om_ks}{2}z^k-\frac{\sin\om_ks}{\om_k}\xi^k+(\frac{\sin\om_k(t-s)+\sin\om_ks}{\sin\om_kt})x^k).$$
Also by $(\ref{H_linverse})$,

$$ H_1^{-1}=\left(%
\begin{array}{ccc}
  D(\frac{-1}{2\vec \om}\cot(\frac{\vec \om t}{2})) & 0 & 0 \\
  0 & 0 & -I \\
  0 & -I & -D(\vec w \cot \vec \om t) \\
\end{array}%
\right).$$
Hence the only non-zero entries of $H_1^{-1}$ are the
ones of the form $h_1^{x^kx^k}, h_1^{z^k \xi ^k}=h_1^{\xi ^k z
^k},$ and $h_1^{\xi ^k \xi ^k}$. Now we let

$$\left\{\begin{array}{c}
  i_{x^kx^k}= \text{the number of times} \,\, h_1^{x^kx^k}\,\, \text{appears in} \,\,
  h_1^{r_1r_2}...h_1^{r_{2j+1}r_{2j+2}} \text{in} \,\, (\ref{P_j+1b_1}),   \\
  i_{z^k \xi ^k}= \text{the number of times} \,\, h_1^{z^k \xi ^k}\,\, \text{appears in} \,\,
  h_1^{r_1r_2}...h_1^{r_{2j+1}r_{2j+2}} \text{in} \,\, (\ref{P_j+1b_1}), \\
  i_{\xi ^k \xi ^k}= \text{the number of times} \,\, h_1^{\xi ^k \xi ^k}\,\, \text{appears in} \,\,
  h_1^{r_1r_2}...h_1^{r_{2j+1}r_{2j+2}} \text{in} \,\, (\ref{P_j+1b_1}). \\
\end{array}\right. $$
By applying these notations to $(\ref{P_j+1b_1})$, $(\ref{b_l})$
we get
\\ \\

$$\! \! \! \!\!\!\!\!\!\!\!\!\!\!\!\!\! P_{j+1}b_1=\frac{i^{-(j+1)}}{2^{j+1}(j+1)!}\sum_{
  \sum_{k=1}^n i_{x^kx^k}+
  i_{z^k \xi ^k}+
  i_{\xi ^k \xi
^k}=j+1}\frac{(j+1)!\,\,2^{\sum_{k=1}^n i_{z^k \xi ^k}
}}{\prod_{k=1}^n i_{x^kx^k}!i_{z^k \xi ^k}!i_{\xi ^k
\xi^k}!}\prod_{k=1}^n \left(\frac{-\cot \frac{\om_k
t}{2}}{2\om_k}\right)^{i_{x^kx^k}}(-1)^{i_{z^k \xi
^k}}\left(-\om_k \cot \om_k t \right)^{i_{\xi ^k \xi ^k}}$$

$$  \qquad \qquad \qquad \times \prod_{k=1}^n \left\{ \left(\frac{\sin\om_k(t-s)+\sin\om_ks}{\sin\om_kt}
\right)^{2i_{x^kx^k}} \left(\frac{\cos \om_ks}{2} \right)^{i_{z^k
\xi ^k}}\left(\frac{-\sin\om_ks}{\om_k} \right)^{i_{z^k \xi
^k}+2i_{\xi ^k \xi^k}} \right \}.
D_{2\al_1,...2\al_n}^{2j+2}W_\hbar,$$
where $\al_k=i_{x^kx^k}+
i_{z^k \xi ^k}+ i_{\xi ^k \xi^k}$, for $k=1,...,n$.
\\

Next we write the above big sum as

$$\!\!\!\!\!\sum_{\sum_{k=1}^n i_{x^kx^k}+i_{z^k \xi ^k}+i_{\xi ^k \xi^k}=j+1}\frac{(j+1)!}
{\prod_{k=1}^n i_{x^kx^k}!i_{z^k \xi ^k}!i_{\xi ^k
\xi^k}!}(\star)= \sum_{\sum\al_k=j+1}\frac{(j+1)!} {\prod \al_k!}
  \sum_{i_{x^kx^k}+i_{z^k \xi ^k}+i_{\xi ^k \xi^k}=\al_k}\prod_{k=1}^n\frac{\al_k!}
{i_{x^kx^k}!i_{z^k \xi ^k}!i_{\xi ^k \xi^k}!}(\star).$$ So the
coefficient of $D_{2\al_1,...2\al_n}^{2j+2}W_\hbar$ in
$P_{j+1}b_1$, equals

$$\frac{i^{-(j+1)}}{2^{j+1}}\frac{(-1)^{j+1}} {(\prod \al_k!)(\prod \om_k^{\al_k})}
\prod_{k=1}^n \left(\half \cot \frac{\om_k
t}{2}\left(\frac{\sin\om_k(t-s)+\sin\om_ks}{\sin\om_kt}\right)^2-\cos
\om_ks\sin\om_k s+\cot\om_kt\sin^2\om_ks \right)^{\al_k}.
$$
Now we observe that the term in the parenthesis simplifies to

$$\half \cot \frac{\om_k
t}{2}\left(\frac{\sin\om_k(t-s)+\sin\om_ks}{\sin\om_kt}\right)^2-\cos
\om_ks\sin\om_k s+\cot\om_kt\sin^2\om_ks=\half \cot \frac{\om_k
t}{2}.$$
So we get
\begin{equation} \label{Pj+1b_1closed}
P_{j+1}b_1=\frac{1}{(2i)^{j+1}} \sum_{|\vec \al|=j+1}
\frac{1}{\vec \al!} \left(\frac{-1}{2\vec \om} \cot \frac{\vec \om
t}{2}\right)^{\vec \al}D_{2 \vec \al}^{2j+2}W_\hbar,
\end{equation}
Finally, by plugging $(x,\vec z,\vec \xi\,\,)=0$ into equation
$(\ref{Pj+1b_1closed})$ and applying it to $(\ref{ajt})$, we get
$(\ref{linearterm})$. This finishes the proof of Theorem
$\ref{Theorem1}$.3.
\\

For future reference let us highlight the equation we just proved

\begin{equation} \label{S_1}
S_1:=\sum_{r_1,...,r_{2j+2}\in \mathcal
A_1}h_1^{r_1r_2}...h_1^{r_{2j+1}r_{2j+2}}\frac{\d
^{2j+2}W}{\d{r_1}...\d{r_{2j+2}}}=(j+1)!\sum_{|\vec \al|=j+1}
\frac{1}{\vec \al!} \left(\frac{-1}{2\vec \om} \cot \frac{\vec \om
t}{2}\right)^{\vec \al}D_{2 \vec \al}^{2j+2}W,
\end{equation}
where $$ W=W(\frac
{\cos\om_ks}{2}z^k-\frac{\sin\om_ks}{\om_k}\xi^k+(\frac{\sin\om_k(t-s)+\sin\om_ks}{\sin\om_kt})x^k).
$$

\subsection{Calculations of $\int_0^t\int_0^{s_1}P_{j+2}b_2(0)$, and the proof of Theorem $\ref{Theorem2}$.}

Throughout this section we assume that $V$ is of the form
($\ref{fg}$). Hence, the only non-zero Taylor coefficients are of
the form $D^{2j+2}_{2\vec \al}V(0)$, or $D^{2j+1}_{2\vec \al+3\vec
e_n}V(0)$, where $\vec e_n=(0,...,0,1)$.

We notice that based on our discussion in the previous section,
the Taylor coefficients of order $2j+1$ appear in
$\int_0^t\int_0^{s_1}P_{j+2}b_2(0)$, and they are of the form
$D^{2j+1}_{\vec \beta}V(0) D^3_{\vec \delta}V(0)$. Therefore we
look for the coefficients of the data

\begin{equation} \label{data} \Big\{D^{2j+1}_{2\vec \al+3\vec
e_n}V(0)D^3_{3\vec e_n}V(0);\quad |\vec \al|=j-1
\Big\}.\end{equation} in the expansion of $a_j(t)$.

\begin{prop} \label{Prop} In the expansion of $a_j(t)$, the coefficient of the data $D^{2j+1}_{2\vec \al+3\vec
e_n}V(0)D^3_{3\vec e_n}V(0)$, $|\vec \al|=j-1$, is

\begin{equation} \label{quadradic}
\frac{c_2(n)}{(2i)^{j+2}}\frac{t}{\vec \al!} \left(\frac{-1}{2\vec
\om} \cot \frac{\vec \om t}{2}\right)^{\vec \al}
\left(\frac{1}{3\om_n^2}(\frac{2\al_n+5}{\al_n+1})(\frac{-1}{2\om_n}
\cot \frac{\om_n t}{2})^2+\frac{1}{9\om_n^4} \right).
\end{equation}
Therefore
\begin{equation} \label{ajt3}
 \!\!\!\!\! \! \!\!\!\!\!\!\!\!\!\!\!\!\!\!\!\!\!\!\!a_j(t)=\frac{c_1(n)}{(2i)^{j+1}}\sum_{|\vec \al|=j+1} \frac{t}{\vec
\al !}(\frac{-1}{2\vec \om}\cot\frac{\vec \om t}{2})^{\vec
\al}D_{2 \vec \al}^{2j+2}V(0) \end{equation}
  $$\qquad \qquad \qquad +\frac{c_2(n)}{(2i)^{j+2}}\sum_{|\vec \al|=j-1}\frac{t}{\vec \al!}
\left(\frac{-1}{2\vec \om} \cot \frac{\vec \om t}{2}\right)^{\vec
\al}
\left(\frac{1}{3\om_n^2}(\frac{2\al_n+5}{\al_n+1})(\frac{-1}{2\om_n}
\cot \frac{\om_n t}{2})^2+\frac{1}{9\om_n^4}
\right)D^{2j+1}_{2\vec \al+3\vec e_n}V(0)D^3_{3\vec e_n}V(0)$$
  $$\!\!\!\!\!\!\!\!\!\!\!\!\!\!+\{\text{a polynomial of Taylor coefficients of order $\leq 2j$}\}.$$
\end{prop}
~~~~~\\

\noindent Before we prove Proposition $\ref{Prop}$, let us show
how to use this proposition to prove Theorem $\ref{Theorem2}$.

 ~~~~~~ \\ \\ \textsc{Proof of Theorem $\ref{Theorem2}$.} First of all, we
prove that for all $\vec \al$, the functions
$$\big(\cot \frac{\vec \om }{2}t\big)^{\vec \al},$$
are linearly independent over $\mathbb C$. To show this we define
$$ \left\{ \begin{array}{ll}
  \vec {\cot}:(0,\pi)^n \longrightarrow \mathbb R^n, \\
   \vec{\cot}(x_1,...,x_n)=(\cot(x_1),...,\cot(x_n)).\\
\end{array} \right.$$
Because $\om_k$ are linearly independent over $\mathbb Q$, the set
$\{(\frac{\om_1}{2}t,...,\frac{\om_n}{2}t)+\pi \mathbb Z^n;\; t\in
\mathbb R\}\cap (0,\pi)^n$ is dense in $(0,\pi)^n$. Since
$\vec{\cot}$ is a homeomorphism and is $\pi$-periodic, we conclude
that the set
$\{(\cot(\frac{\om_1}{2}t),...,\cot(\frac{\om_n}{2}t);\, t\in
\mathbb R\}$ is dense in $\mathbb R^n$. Now assume
$$ \sum_{\vec \al}c_{\vec \al}\big(\cot \vec
{\frac{\om }{2}}t\big)^{\vec \al}=0.$$ Since
$\{(\cot(\frac{\om_1}{2}t),...,\cot(\frac{\om_n}{2}t);\, t\in
\mathbb R\}$ is dense in $\mathbb R^n$, we get
$$\sum_{\vec \al}c_{\vec \al}\vec X^{\vec \al}=0,$$
for every $\vec X=(X_1,...,X_n)\in \mathbb R^n$. But the monomials
$\vec X^{\vec \al}$ are linearly independent over $\mathbb C$. So
$c_{\vec \al}=0$.
\\

Next we argue inductively to recover the Taylor coefficients of
$V$ from the wave invariants. Since
$$a_0(t)=\prod_{k=1}^n \frac{1}{2i\sin \frac{\om_k t}{2}},$$
we can recover $\prod_{k=1}^n \sin \frac{\om_k t}{2}$, and
therefore we can recover $\{\om_k\}$ up to a permutation. This can
be seen by Taylor expanding $\prod_{k=1}^n \sin \frac{\om_k
t}{2}$. We fix this permutation and we move on to recover the
third order Taylor coefficient $D^3_{3\vec e_n}V(0)$. This term
appears first in $a_1(t)$. By Proposition \ref{Prop}, we have

$$a_1(t)=c_1(n)\frac{t}{(2i)^{2}}\sum_{|\vec \al|=2} \frac{1}{\vec
\al !}(\frac{-1}{2\vec \om}\cot\frac{\vec \om}{2}t)^{\vec \al}D_{2
\vec \al}^{4}V(0)$$
  $$\qquad \qquad \qquad + \,c_2(n)\frac{t}{(2i)^{3}}
\left(\frac{5}{3\om_n^2}(\frac{-1}{2\om_n} \cot \frac{\om_n
t}{2})^2+\frac{1}{9\om_n^4} \right)\Big(D^{3}_{3\vec
e_n}V(0)\Big)^2
$$
  $$\!\!\!\!\!\!\!\!\!\!\!\!\!\!+\{\text{a rational function of}\; \om_k\}.$$
Now since the functions $\{(\cot\frac{\vec \om}{2}t)^{\vec
\al}\}_{|\vec \al|=2}$ and
$\left(\frac{5}{3\om_n^2}(\frac{-1}{2\om_n} \cot \frac{\om_n
t}{2})^2+\frac{1}{9\om_n^4} \right)$ are linearly independent over
$\mathbb C$, we can therefore recover the data $\{D_{2 \vec
\al}^{4}V(0)\}_{|\vec \al|=2}$ and $\{D^{3}_{3\vec e_n}V(0)^2\}$
from $a_1(t)$. So we have determined the third order term $
D^{3}_{3\vec e_n}V(0)$ up to a minus sign from the first invariant
$a_1(t)$. This choice of minus sign corresponds to a reflection.
We fix this reflection and we move on to determine the higher
order Taylor coefficients inductively.

Next we assume $ D^{3}_{3\vec e_n}V(0)\neq 0$ and that we know all
the Taylor coefficients $D^m_{\vec \beta}V(0)$ with $m\leq2j$. We
wish to determine the data $\{D^{2j+1}_{2 \vec \al +3 \vec
e_n}V(0)\}_{|\vec \al|=j-1}$ and $\{D^{2j+2}_{2\vec
\al}V(0)\}_{|\vec \al|=j+1}$, from the wave invariant $a_j(t)$. At
this point we use Proposition \ref{Prop}, and to finish the proof
of Theorem \ref{Theorem1} we have to show that the set of
functions

$$\Big\{(\cot\frac{\vec \om}{2}t)^{\vec \al}; \; \quad |\vec \al|=j+1\Big\}\cup
\Big\{(\cot \frac{\vec \om t}{2})^{\vec \al}
\left(\frac{1}{3\om_n^2}(\frac{2\al_n+5}{\al_n+1})(\frac{-1}{2\om_n}
\cot \frac{\om_n t}{2})^2+\frac{1}{9\om_n^4} \right); \; \quad
|\vec \al |=j-1 \Big\},$$ are linearly independent over $\mathbb
C$. But this is clear from our discussion at the beginning of
proof.
\\

\noindent \textsc{Proof of Proposition $\ref{Prop}$.} As we
mentioned at the beginning of Section $2.7$, the data
$D^{2j+1}_{2\vec \al+3\vec e_n}V(0)D^3_{3\vec e_n}V(0)$, $|\vec
\al|=j-1$, appears first in $a_j(t)$ and it is a part of the term
$\int_0^t\int_0^{s_1}P_{j+2}b_2(0)$. So let us calculate those
terms in the expansion of $P_{j+2}b_2(0)$ which contain
$D^{2j+1}_{2\vec \al+3\vec e_n}V(0)D^3_{3\vec e_n}V(0)$. By
(\ref{b_l}), since here $l=2$, we have

\begin{equation} \label{b_2} \begin{array}{lll}
  b_2(s_1,s_2,x,x,z_1,z_2,\xi_1,\xi_2)=W_1\,W_2, \qquad \text{where,} \\ \\
 W_1= W(\frac
{\cos\om_ks_1}{2}({z_1^k+z_2^k})-\frac{\sin\om_ks_1}{\om_k}\xi_1^k+(\frac{\sin\om_k(t-s_1)+\sin\om_ks_1}{\sin\om_kt})x^k),\\ \\
   W_2=W(\frac
{\cos\om_ks_2}{2}z_2^k-\frac{\sin\om_ks_2}{\om_k}\xi_2^k+(\frac{\sin\om_k(t-s_2)+\sin\om_ks_2}{\sin\om_kt})x^k).\\
\end{array}\end{equation}
Also from (\ref{H_linverse}) we have

\begin{equation} \label{H_2inverse}
H_2^{-1}=
\left(%
\begin{array}{cccccc}
  D(\frac{-1}{2\vec \om}\cot(\frac{\vec \om t}{2})) & 0 & 0 & 0 &0
  \\ \\
  0 & 0 & 0& -I & -I \\ \\
  0 & 0 & 0& 0 & -I \\ \\
  0 & -I & 0 & -D(\vec \om\cot(\vec \om t)) & -D(\vec \om\cot(\vec \om
  t))\\ \\
  0 & -I & -I & -D(\vec \om\cot(\vec \om t)) & -D(\vec \om\cot(\vec \om
  t))\\ \\
\end{array}%
\right)_{5n\times 5n}. \\
\end{equation}
By (\ref{P_jb_l}) and (\ref{b_2}), $P_{j+2}b_2(0)=
\frac{i^{-j}}{2^j.j!}S_2$, where $S_2$ is the following sum

\begin{equation} \label{S_2} S_2= \sum_{r_1,...,r_{2j+4}\in\mathcal
A_2}h_2^{r_1r_2}...h_2^{r_{2j+3}r_{2j+4}}(W_1\,W_2)_{{r_1}...{r_{2j+4}}}(0),
\end{equation}
where $\mathcal A_2=\{x^k, z_1^k, z_2^k, \xi_1^k,
\xi_2^k\}_{k=1}^n$ and for every $r,r' \in \mathcal A_2$,
$h_2^{rr'}$ is the $(r,r')$-entry of the matrix $H_2^{-1}$ in
(\ref{H_2inverse}). We would like to separate out those terms in
$S_2$ which include $D^{2j+1}_{2\vec \al+3\vec e_n}V(0)D^3_{3\vec
e_n}V(0)$. To do this, from the total number $2j+4$ derivatives
that we want to apply to $W_1W_2$, we have to put $3$ of them on
$W_1$ \, (or $W_2$) and put $2j+1$ of them on $W_2$ (or $W_1$
respectively). These combinations fit into one of the following
two different forms
\begin{equation} \label{S_2^1} S_2^1=
\sum_{r_1,...,r_{2j+4}\in\mathcal
A_2}h_2^{r_1r_2}...h_2^{r_{2j+3}r_{2j+4}}(W_1)_{r_1r_2r_3}\,(W_2)_{r_4,r_5...{r_{2j+4}}}(0).
\end{equation}
There are $2(j+1)(j+2)$ terms of this form in the expansion of
$S_2$.
\begin{equation} \label{S_2^2} S_2^2=
\sum_{r_1,...,r_{2j+4}\in\mathcal
A_2}h_2^{r_1r_2}...h_2^{r_{2j+3}r_{2j+4}}(W_1)_{r_1r_3r_5}\,(W_2)_{r_2,r_4,r_6,r_7...{r_{2j+4}}}(0).
\end{equation}
There are $2^3\left(\begin{array}{c}
  j+2 \\
  3 \\
\end{array}\right)$ terms of this form in the expansion of
$S_2$.
\\

\noindent Now, we calculate the sums $S_2^1$ and $S_2^2$.
\\

\subsection{Calculation of $S_2^1$.} We rewrite $S_2^1$ as

$$S_2^1=
\sum_{r_1,...,r_{4}}h_2^{r_1r_2}h_2^{r_3r_4}\Big(\sum_{r_5,...,r_{2j+4}}h_2^{r_5r_6}...h_2^{r_{2j+3}r_{2j+4}}
(W_2)_{r_5...{r_{2j+4}}}\Big)_{r_4}(W_1)_{r_1r_2r_3}(0).
$$
Then from the definition of $W_2$ in (\ref{b_2}) and also from
(\ref{H_2inverse}) it is clear that we can apply (\ref{S_1}) to
the sum in the big parenthesis above. Hence we get

\begin{equation} \label{S_2^1new} S_2^1=(j+1)!\sum_{|\vec \al|=j} \frac{1}{\vec \al!}
\left(\frac{-1}{2\vec \om} \cot \frac{\vec \om t}{2}\right)^{\vec
\al}\Big(\sum_{r_1,...,r_{4}}h_2^{r_1r_2}h_2^{r_3r_4}
\big(D^{2j}_{2\vec
\al}W_2\big)_{r_4}(W_1)_{r_1r_2r_3}\Big)(0).\end{equation} This
reduces the calculation of $S_2^1$ to calculating the small sum

$$A_2^1=\sum_{r_1,...,r_{4}}h_2^{r_1r_2}h_2^{r_3r_4}
(\tilde W_2)_{r_4}(W_1)_{r_1r_2r_3}(0) ,\qquad (\tilde
W_2=D^{2j}_{2\vec \al}W_2).$$ Computation of the sum $A_2^1$ is
straight forward and we omit writing the details of this
computation. Using Maple, we obtain

$$ \int_0^t\int_0^{s_1} A_2^1 \,ds_2\,dt=-\frac{t}{2\om_n^2}(\frac{-1}{2\om_n}
\cot \frac{\om_n t}{2})\big(D^1_{\vec e_n}\tilde W_2 \, D^3_{3
\vec e_n} W_1\big)(0).$$ If we plug this into (\ref{S_2^1new}),
after a change of variable $\al_n \rightarrow \al_n+1$ in indices,
we get

\begin{equation} \label{S_2^1nnew}
\int_0^t\int_0^{s_1} S_2^1
\,ds_2\,dt=\frac{(j+1)!}{\al_n+1}\sum_{|\vec \al|=j-1}
\frac{1}{\vec \al!} \left(\frac{-1}{2\vec \om} \cot \frac{\vec \om
t}{2}\right)^{\vec \al}(-\frac{t}{2\om_n^2})(\frac{-1}{2\om_n}
\cot \frac{\om_n t}{2})^2 D^{2j+1}_{2\vec \al+3\vec e_n}V(0) \,
D^3_{3 \vec e_n} V(0).
\end{equation}

\subsection{Calculation of $S_2^2$.}We rewrite $S_2^2$ as
$$S_2^2=\sum_{r_1,...,r_{6}}h_2^{r_1r_2}h_2^{r_3r_4}h_2^{r_5r_6}\Big(\sum_{r_7,...,r_{2j+4}}h_2^{r_7r_8}...h_2^{r_{2j+3}r_{2j+4}}
(W_2)_{r_7...{r_{2j+4}}}\Big)_{r_2,r_4,r_6}(W_1)_{r_1r_3r_5}(0).
$$
Again from (\ref{H_2inverse}) it is clear that we can apply
(\ref{S_1}) to the sum in the big parenthesis above. So

\begin{equation} \label{S_2^2new} S_2^2=(j+1)!\sum_{|\vec \al|=j-1} \frac{1}{\vec \al!}
\left(\frac{-1}{2\vec \om} \cot \frac{\vec \om t}{2}\right)^{\vec
\al}\Big(\sum_{r_1,...,r_{6}}h_2^{r_1r_2}h_2^{r_3r_4}
\big(D^{2j-2}_{2\vec
\al}W_2\big)_{r_2,r_4,r_6}(W_1)_{r_1r_3r_5}\Big)(0).\end{equation}
So we need to compute

$$A_2^2=\sum_{r_1,...,r_{6}}h_2^{r_1r_2}h_2^{r_3r_4}h_2^{r_5r_6}
(\tilde W_2)_{r_2,r_4,r_6}(W_1)_{r_1r_3r_5}(0) ,\qquad (\tilde
W_2=D^{2j-2}_{2\vec \al}W_2).$$ Using Maple
$$ \int_0^t\int_0^{s_1} A_2^2 \,ds_2\,dt=\big(-\frac{t}{2\om_n^2}\big(\frac{-1}{2\om_n}
\cot \frac{\om_n
t}{2}\big)^2-\frac{t}{12\om_n^4}\big)\big(D^3_{\vec 3e_n}\tilde
W_2 \, D^3_{3 \vec e_n} W_1\big)(0).$$ If we plug this into
(\ref{S_2^2new}) we get

\begin{equation} \label{S_2^2nnew}
\int_0^t\int_0^{s_1} S_2^2 \,ds_2\,dt=(j+1)!\sum_{|\vec \al|=j-1}
\frac{1}{\vec \al!} \left(\frac{-1}{2\vec \om} \cot \frac{\vec \om
t}{2}\right)^{\vec
\al}\big(-\frac{t}{2\om_n^2}\big(\frac{-1}{2\om_n} \cot
\frac{\om_n t}{2}\big)^2-\frac{t}{12\om_n^4}\big)D^{2j+1}_{2\vec
\al+3\vec e_n}V(0) \, D^3_{3 \vec e_n} V(0).
\end{equation}
\\
\noindent We note that the part of the expansion of
$\int_0^t\int_0^{s_1}P_{j+2}b_2(0)$ which contains the data
$D^{2j+1}_{2\vec \al+3\vec e_n}V(0) \, D^3_{3 \vec e_n} V(0)$,
equals

$$ \frac{i^{-j}}{2^j.j!}\Big( 2(j+2)(j+1)\int_0^t\int_0^{s_1}S_2^1+2^3\big(\begin{array}{c}
  j+2 \\
  3 \\
\end{array}\big) \int_0^t\int_0^{s_1}S_2^2\Big).$$ Finally, by applying equations (\ref{S_2^1nnew}) and
(\ref{S_2^2nnew}) to this we obtain (\ref{ajt3}).

\section{Appendix A}

In this appendix we prove Lemma \ref{criticallem}.
\\

\textsc{Proof.} First of all we would like to change the function
$\Theta$ slightly by rescaling it. We choose $0<\tau <2\ep$ so
that $\hbar^{1-\tau}=o(\hbar^{1-2\ep})$. Then we define

\begin{equation} \label{Thetah}\Theta_\hbar(x):=\Theta(\frac{x}{\hbar^{1-\tau}}).
\end{equation}
Thus $\Theta_\hbar\in C^\infty_0([0,\infty))$ is supported in the
interval $I_\hbar=[0, \hbar^{1-\tau}\delta]$. In Appendix B, using
min-max principle we show that

$$Tr(\Theta(\hat H)e^{\frac{-it}{\hbar} \hat
H})=Tr(\Theta_\hbar(\hat H)e^{\frac{-it}{\hbar} \hat H})+O(\hbar
^\infty)=Tr(e^{\frac{-it}{\hbar} \hat H})+O(\hbar ^\infty).$$

Hence to prove the lemma it is enough to show

$$ Tr(\Theta(\hat P)e^{\frac{-it}{\hbar} \hat P})=Tr(\Theta_\hbar(\hat H)e^{\frac{-it}{\hbar} \hat
H})+O(\hbar ^\infty).$$

To prove this identity we use the WKB construction of the kernel
of the operators $\Theta(\hat P)e^{\frac{-it}{\hbar} \hat P}$ and
$\Theta_\hbar(\hat H)e^{\frac{-it}{\hbar} \hat H}$ and make a
compression between them.

\subsection{WKB construction for $\Theta(\hat P)e^{\frac{-it}{\hbar} \hat
P}$}

In \cite{DSj}, chapter $10$, a WKB construction is made for
$\Theta(\hat P)e^{\frac{-it}{\hbar} \hat P}$ for symbols $P$ in
the symbol class $S_0^0(1)$ which are independent of $\hbar$ or of
the form $P(x,\xi,\hbar) \sim P_0(x,\xi)+\hbar P_1(x,\xi)+...,$
where $P_j\in S_0^0$ are independent of $\hbar$ (but not for
symbols $H=H(x,\xi,\hbar)\in S^0_{\delta_0}$).

It is shown that we can approximate $\Theta(\hat
P)e^{\frac{-it}{\hbar} \hat P}$ for small time $t$, say
$t\in(-t_0,t_0)$, by a fourier integral operator of the form

$$ U_{P}(t)u(x)=(2\pi\hbar)^{-n}\int\int
e^{i(\phi_P(t,x,\eta)-y.\eta)/\hbar}
b_P(t,x,y,\eta,\hbar)u(y)dyd\eta,$$ where $b_P\in C^\infty
((-t_0,t_0);\;S(1))$ have uniformly compact support in
$(x,y,\eta)$, and $\phi_P$ is real, smooth and is defined near the
support of $b_P$. The functions $\phi_P$ and $b_P$ are found in
such a way that for all $t\in (-t_0,t_0)$

$$ ||\Theta(\hat P)e^{\frac{-it}{\hbar} \hat
P}-U_{P}(t)||_{tr}=O(\hbar^{\infty}).$$

Let us briefly review this construction, made in \cite{DSj}. First
of all, in Chapter 8, Theorem 8.7, it is proved that for every
symbol $P\in S_0^0(1)$, we have $\Theta(\hat
P)=Op^w(a_{P}(x,\xi,\hbar))$ for some $a_{P}(x,\xi,\hbar)\in
S_0^0(1)$, where here $\hat P$ and $Op^w(a_{P}(x,\xi,\hbar))$ are
respectively the Weyl quantization of $P$ and
$a_{P}(x,\xi,\hbar)$. It is also shown that $a_{P} \sim
a_{P,0}(x,\xi)+ha_{P,1}(x,\xi)+...$ for some $a_{P,j}(x,\xi)\in
S_0^0(1)$. The idea of proof is as follows. In Theorem $8.1$ of
\cite{DSj} it is shown that if $\Theta \in C_0^{\infty}(\mathbb
R)$, and if $\tilde \Theta\in C^1_0(\mathbb C)$ is an almost
analytic extension of $\Theta$ (i.e. $\bar \d \tilde
\Theta(z)=O(|\Im z|^\infty)),$ then

\begin{equation} \label{Theta}
\Theta(\hat P)=\frac{-1}{\pi} \int_{\mathbb C} \frac{\dbar \tilde
\Theta (z)}{z-\hat P}L(dz).
\end{equation}
Then it is verified that for some symbol $r(x,\xi,z;\hbar)$, we
have $(z-\hat P)^{-1}={Op}^w(r(x,\xi,z;\hbar))$. By symbolic
calculus, one can find a formal asymptotic expansion of the form

$$ r(x,\xi,z;\hbar)\sim \frac{1}{z- P}+ \hbar\frac{q_1(x,\xi,z)}{(z-
P)^3}+\hbar^2\frac{q_2(x,\xi,z)}{(z- P)^5}+...,$$
by formally
solving $Op^w(r(x,\xi,z;\hbar))\sharp_\hbar (z-\hat P)=(z-\hat
P)\sharp_\hbar Op^w(r(x,\xi,z;\hbar))=1$. We can see that
$q_j(x,\xi,z)$ are polynomials in $z$ with smooth coefficients.
Finally it is shown that $\Theta(\hat
P)=Op^w(a_{P}(x,\xi,\hbar))$, where $a_P \in S_0^0$ is given by

$$a_{P}(x,\xi,\hbar)=\frac{-1}{\pi}\int_{\mathbb C} \dbar \tilde
\Theta (z)r(x,\xi,z;\hbar)L(dz).$$ By the above asymptotic
expansion for $r(x,\xi,z;\hbar)$ one obtains an asymptotic
$a_{P}\sim a_{P,0}+\hbar a_{P,1}+...$, where

\begin{equation} \label{aPj}
a_{P,j}=\frac{-1}{\pi}\int_{\mathbb C} \dbar \tilde \Theta
(z)\frac{q_j(x,\xi,z)}{(z-P)^{2j+1}}L(dz)=\frac{1}{(2j)!}\d_t^{2j}(q_j(x,\xi,t)
\Theta(t))|_{t=P(x,\eta)}.
\end{equation}
Then, again in Chapter 10 of \cite{DSj}, it is shown that
$\phi_P(t,x,\eta)$ and $b_P(t,x,y,\eta,\hbar)$ satisfy
\begin{equation} \label{eikonalP}\d_t\phi_P(t,x,\eta)+P(x,\d_x\phi_P(t,x,\eta))=0, \qquad
\phi_P|_{t=0}=x.\eta, \end{equation}
$$b_P\sim b_{P,0}+\hbar b_{P,1}+.\, .\, .\, , \quad \qquad
b_{P,j}=b_{P,j}(t,x,y,\eta) \in C^\infty((-t_0,t_0); \;
S^0_{0}(1)),$$ where
\begin{equation} \label{transportP}\left\{ \begin{array}{lll}
  \d_t b_{P,j}+ \big<\d_x \phi_P,\d_xb_{P,j}\big>+\half\Delta_x
  \phi_P
  .\,b_{P,j}=-\half\Delta_x
  b_{P,j-1}, \qquad j\geq 0,\qquad (b_{P,-1}=0),
  \\ \\
 b_{P,j}|_{t=0}=\psi(x,\eta)a_{P,j}(\frac{x+y}{2},\eta)\psi(y,\eta). \\
\end{array} \right. \end{equation}
In (\ref{transportP}), $a_{P,j}$ is given by (\ref{aPj}) and
$\psi(x,\eta)$ is any $C_0^\infty$ function which equals $1$ in a
neighborhood of $\overline{P^{-1}(I)}$ where $I=[0,\delta]$ is, as
before, the range of our low-lying eigenvalues and where $\Theta$
is supported.
\\

\noindent There exists a similar construction for
$\Theta_\hbar(\hat H)e^{\frac{-it}{\hbar} \hat H}$, except here
$H\in S^0_{\delta_0}$.

\subsection{WKB construction for $\Theta_\hbar(\hat
H)e^{\frac{-it}{\hbar} \hat H}$}

Since in (\ref{PH}), $H=H(x,\xi,\hbar)\in S^0_{\delta_0}$, with
$\delta_0=\half-\ep$, we can not simply use the construction in
\cite{DSj} mentioned above. Here in two lemmas we show that the
same construction works for the operator $\Theta_\hbar(\hat
H)e^{\frac{-it}{\hbar} \hat H}$. We will closely follow the proofs
in \cite{DSj}.

\begin{lem}\label{ThetaH}
\quad

\begin{itemize}

\item[1)] Let $\Theta_\hbar$ be given by (\ref{Thetah}) and $H \in
S^0_{\delta_0}$ by (\ref{PH}). Then for some $a_H\in
S^0_{\delta_0}$ we have $\Theta(\hat H)=Op^w(a_{H}(x,\xi,\hbar))$.
Moreover $a_H(x,\xi,\hbar)\sim a_{H,0}(x,\xi,\hbar)+\hbar
a_{H,1}(x,\xi,\hbar)+...,$ where $a_{H,j}(x,\xi,\hbar)\in
S^0_{\delta_0}$ is given by

\begin{equation}\label{aHj}
a_{H,j}=\frac{-1}{\pi}\int_{\mathbb C} \dbar \tilde \Theta
(z)\frac{q_{H,j}(x,\xi,z,\hbar)}{(z-H)^{2j+1}}L(dz)=\frac{1}{(2j)!}\d_t^{2j}(q_{H,j}(x,\xi,t,\hbar)
\Theta_\hbar(t))|_{t=H(x,\xi,\hbar)}.
\end{equation}

\item[2)] Choose $c$ such that $0<c<\min\{1,\om_k^2\}_{k=1}^n \leq
\max\{1,\om_k^2\}_{k=1}^n <\frac{1}{c}$. Let $\psi_\hbar(x,\eta)$
be a function in $C^\infty_0(\mathbb R^{2n}) \cap
S^0_{\delta_0}(\mathbb R^{2n})$ which is supported in the ball
$\{x^2+\eta^2<4c^{-1}\hbar^{1-\tau} \delta\}$ and equals $1$ in a
neighborhood of $\overline{H^{-1}(I)}$, where
$I_\hbar=[0,\hbar^{1-\tau}\delta]$\,($I_\hbar$ is where
$\Theta_\hbar$ is supported). Then

\begin{equation} \label{microlocalTheta}\Theta_\hbar(\hat H)u(x)=(2\pi \hbar)^{-n} \int\int
e^{i(x-y).\eta
/\hbar}\psi_\hbar(x,\eta)a_{H}(\frac{x+y}{2},\eta,\hbar)\psi_\hbar(y,\eta)u(y)dyd\eta
+K(\hbar)u(x), \end{equation}

\noindent where $||K(\hbar)||_{tr}=O(\hbar ^{\infty})$.
\end{itemize}
\end{lem}

\textsc{Proof of Lemma \ref{ThetaH}}: Since $H\in S^0_{\delta_0}$
and $\delta_0=\half-\ep<\half$, the symbolic calculus mentioned in
the last section can be followed similarly to prove Lemma
\ref{ThetaH}.1. It is also easy to check that in (\ref{aHj}),
$a_{H,j}\in S^0_{\delta_0}$. The second part of the Lemma is
stated in \cite{DSj} , equation $10.1$, for the case $P\in S_0^0$.
The same argument works for $H\in S_{\delta_0}^0$, precisely
because the factor $\hbar^N$ on the right hand side of the
inequality in Proposition $9.5$ of \cite{DSj} changes to
$\hbar^{N-\delta_0 \al}$. Thus the discussion on pages $115-116$
still follows.
\\

\begin{lem} \label{UH}For every $t$ in some small interval $(-t_0,t_0)$,
there exist functions $\phi_H(t,x,\eta,\hbar)$ and
$b_H(t,x,y,\eta,\hbar)$ such that the operator $U_H(t)$ defined by

\begin{equation}\label{UH1}
U_{H}(t)u(x)=(2\pi\hbar)^{-n}\int\int
e^{i(\phi_H(t,x,\eta,\hbar)-y.\eta)/\hbar}
b_H(t,x,y,\eta,\hbar)u(y)dyd\eta, \end{equation}

satisfies $$||\Theta_\hbar(\hat H)e^{\frac{-it}{\hbar} \hat
H}-U_H(t)||_{tr}=O(\hbar^\infty).$$
\\

Moreover, we can choose $\phi_H$ and $b_H$ such that

\begin{itemize}

\item[1)] $\phi_H$ satisfies the eikonal equation
\begin{equation}\label{eikonalH}
\d_t\phi_H(t,x,\eta,\hbar)+H(x,\d_x\phi_H(t,x,\eta,\hbar))=0,
\qquad \phi_H|_{t=0}=x.\eta.
\end{equation}

\noindent This equation can be solved in $(-t_0,t_0) \times
\{x^2+\eta^2<C\hbar^{1-\tau}\delta\}$ where $C$ is an arbitrary
constant. In fact $\phi_H$ is independent of $\hbar$ in this
domain. (Only the domain of $\phi_H$ depends on $\hbar$. See
(\ref{phi_Hphi_P}).)
\\

\item[2)] For all $t\in(-t_0,t_0)$, we have $b_H(t,x,y,\eta,\hbar)
\in S^0_{\delta_0}$ with supp\,$b_H \subset \{x^2+\eta^2,
y^2+\eta^2 <C_1\hbar^{1-\tau}\delta\}$ for some constant $C_1$.
Also $b_H$ has an asymptotic expansion of the form

\begin{equation} \label{bH} b_H\sim b_{H,0}+\hbar b_{H,1}+.\, .\, .\, , \quad \qquad
b_{H,j}=b_{H,j}(t,x,y,\eta,\hbar) \in C^\infty((-t_0,t_0); \;
S^0_{\delta_0}(1)), \end{equation}

\noindent and the functions $b_{H,j}$ satisfy the transport
equations

\begin{equation}\label{transportH}\left\{ \begin{array}{lll}
  \d_t b_{H,j}+ \big<\d_x \phi_H,\d_xb_{H,j}\big>+\half\Delta_x \phi_H
  .\,b_{H,j}=-\half\Delta_x
  b_{H,j-1}, \qquad j\geq 0,\qquad (b_{H,-1}=0),
  \\ \\
 b_{H,j}|_{t=0}=\psi_\hbar(x,\eta)a_{H,j}(\frac{x+y}{2},\eta,\hbar)\psi_\hbar(y,\eta), \\
\end{array} \right. \end{equation}

\noindent where in (\ref{transportH}) we let $\psi_\hbar(x,\eta)$
be a function in $C^\infty_0(\mathbb R^{2n}) \cap
S^0_{\delta_0}(\mathbb R^{2n})$ which is supported in the ball
$\{x^2+\eta^2<4c\hbar^{1-\tau} \delta\}$ and equals $1$ in a
neighborhood of $\overline{H^{-1}(I_\hbar)}$, where
$I_\hbar=[0,\hbar^{1-\tau}\delta]$. Here $c$ is defined in Lemma
\ref{ThetaH}.2. Also in (\ref{transportH}), the functions
$a_{H,j}$ are defined by (\ref{aHj}).

\item[4)]\begin{equation} \label{phi_Hphi_P} \qquad \qquad
\phi_H(t,x,\eta,\hbar)=\phi_P(t,x,\eta)\qquad \text{on}\qquad
\{x^2+\eta^2, y^2+\eta^2 <C_1\hbar^{1-\tau}\delta\}\supset
\text{supp}\,(b_H(x,y,\eta,\hbar)).
\end{equation}

\item[5)] \begin{equation}\label{b_Hb_P}
b_{H,j}(t,x,y,\eta,\hbar)=b_{P,j}(t,x,y,\eta) \qquad
\text{on}\qquad \big\{x^2+\eta^2,\, y^2+\eta^2 < c
\hbar^{1-\tau}\delta \big\}.
\end{equation}

\end{itemize}
\end{lem}

\textsc{Proof of Lemma \ref{UH}:} First of all we assume $U_H(t)$
is given by (\ref{UH1}) and we try to solve the equation

$$\left\{\begin{array}{lll}
  ||(\frac{\hbar}{i}\d_t+\hat H)U_H(t)||_{tr}=O(\hbar^\infty), \\ \\
  U_{H}(0)=\Theta_\hbar(\hat H). \\
\end{array} \right.$$

for $\phi_H$ and $b_H$, for small time $t$.

Using (\ref{microlocalTheta}), this leads us to

$$ \left\{\begin{array}{lll}
  e^{-i\phi_H/\hbar}(\frac{\hbar}{i}\d_t+\hat H)(e^{i\phi_H/\hbar}b_H)\in C^\infty((-t_0,t_0);\; S^{-\infty}_{\delta_0}(1)), \\
  \\
b|_{t=0}= \psi_\hbar(x,\eta)a_{H}(\frac{x+y}{2},\eta,\hbar)\psi_\hbar(y,\eta). \\
\end{array} \right.$$

We choose the phase function $\phi_H=\phi_H(t,x,\eta,\hbar)$ to
satisfy the eikonal equation (\ref{eikonalH}). This equation can
be solved in a neighborhood of the support of $b_H$, for small
time $t\in(-t_0,t_0)$ with $t_0$ independent of $\hbar$. Let us
explain how to solve this equation. We let $(x(t,z,\eta;\hbar),
\xi(t,z,\eta;\hbar))$ be the solution to the Hamilton equation

\begin{equation}\label{Hamilton}\left\{\begin{array}{lll}
  \d_tx=\d_\xi H(x,\xi,\hbar)=\xi, \qquad \qquad \quad \qquad \qquad x(0,z,\eta;\hbar)=z
\\ \\
 \d_t \xi=-\d_x H(x,\xi,\hbar)=-\d_x V_\hbar(x), \qquad \qquad \xi(0,z,\eta;\hbar)=\eta \\
\end{array} \right. . \end{equation}
We can show that (see section 4 of \cite{Ch}) there exists $t_0$
independent of $\hbar$ such that for all $|t|\leq t_0$ we have

\begin{equation} \label{dxdz} \left\{
\begin{array}{ll}
  |\d_z x(t,z,\eta;\hbar)-I|\leq \half,\qquad |\d_\eta x(t,z,\eta;\hbar)|\leq \half
  \\ \\
   |\d_z \xi(t,z,\eta;\hbar)|\leq \half, \qquad |\d_\eta \xi(t,z,\eta;\hbar)-I|\leq \half\\
\end{array}\right. .
\end{equation}
We can choose $t_0$ independent of $\hbar$, precisely because in
equation 4.4 of \cite{Ch} we have a uniform bound in $\hbar$ for
Hess$(V_\hbar(x))$. Now, we define
$$\lambda: (z,\eta) \longmapsto (x(t,z,\eta;\hbar),\eta).$$
It is easy to see that $\la(0,0)=(0,0)$. This is because if
$(z,\eta)=(0,0)$ then $H(x,\xi)=H(z,\eta)=0$. By (\ref{PH}) and
(\ref{Hypothesis}), and $W(x)=O(|x|^3)$, we can see that
$H(x,\xi)=0$ implies $(x(t,0,0;\hbar),\xi(t,0,0;\hbar))=(0,0)$. On
the other hand from (\ref{dxdz}) we have $\half<|\d_z
x(t,z,\eta;\hbar)|<\frac{3}{2}$. Therefore $\lambda$ is invertible
in a neighborhood of origin. We define the inverse function by

$$\lambda^{-1}(x,\eta)=(z(t,x,\eta;\hbar),\eta),$$
which is defined in a neighborhood of $(x,\eta)=(0,0)$. Then we
have
\begin{equation} \label{solveeikonal}
\phi_H(t,x,\eta,\hbar)=z(t,x,\eta;\hbar).\eta\,+\, \int_0^t \half
|\xi(s,z(t,x,\eta;\hbar),\eta;\hbar)|^2-V_\hbar(x(s,z(t,x,\eta;\hbar),\eta;\hbar))ds,
\end{equation}
A similar formula holds for $\phi_P$ except in (\ref{Hamilton})
$H$ should be replaced by $P$ and in (\ref{solveeikonal})
$V_\hbar$ by $V$. It is known that the eikonal equation for
$\phi_P$ can be solved near supp$b_P$, for small time
$t\in(-t_0,t_0)$ \, (Of course $t_0$ is independent of $\hbar$).
Now, we want to show that
\begin{equation} \label{phiphi}
\phi_H(t,x,\eta,\hbar)=\phi_P(t,x,\eta)\qquad \text{in} \qquad
(-t_0,t_0) \times \{x^2+\eta^2<C\hbar^{1-\tau}\delta\}.
\end{equation}
Let $(x,\eta)$ be in $\{x^2+\eta^2<C\hbar^{1-\tau}\delta\}$.
First, we show that
$|z(t,x,\eta;\hbar)|<8C^{\half}\hbar^{\frac{1-\tau}{2}}\delta^{\half}$.
Because $z(t,0,0;\hbar)=0$, by Fundamental Theorem of Calculus we
have
$$|z(t,x,\eta;\hbar)|\leq\big(|x|+|\eta|\big)\text{sup}\{(|\d_x|+|\d_\eta|)(z(t,x,\eta;\hbar))\}.$$
From $x(t,z(t,x,\eta;\hbar),\eta;\hbar)=x$, we get
$$\d_\eta z= -(\d_z x)^{-1}\d_\eta x .$$ Thus by (\ref{dxdz}), $|\d_x
z|+|\d_\eta z|\leq 4$. Hence $|z(t,x,\eta;\hbar)|<
4(|x|+|\eta|)<8C^{\half}\hbar^{\frac{1-\tau}{2}}\delta^{\half}$.
This implies that for all $|t|\leq t_0$,
$(x(s,z(t,x,\eta;\hbar),\eta;\hbar),
\xi(s,z(t,x,\eta;\hbar),\eta;\hbar))$ will stay in a ball of
radius $O(\hbar^{1-\tau})$ centered at the origin. On the other
hand, by definition (\ref{PH}), $P$ and $H$ agree in the ball
$\{x^2+\eta^2<\frac{1}{4}\hbar^{1-2\ep}\}$ and $\tau<2\ep$. So for
all $t,s\in (-t_0,t_0)$ and $(x,\eta)\in
\{x^2+\eta^2<C\hbar^{1-\tau}\delta\}$ we have

\begin{equation} \label{zPzH}\begin{array}{c}
  z_P(t,x,\eta)=z_(t,x,\eta;\hbar), \\ \\
  x_P(s,z_P(t,x,\eta),\eta)=x(s,z(t,x,\eta;\hbar),\eta;\hbar),\\ \\
  \xi_P(s,z_P(t,x,\eta),\eta)=\xi(s,z(t,x,\eta;\hbar),\eta;\hbar),\\
\end{array} \end{equation}
where $z_P(t,x,\eta)$, $x_P(s,z_P(t,x,\eta),\eta)$ and
$\xi_P(s,z_P(t,x,\eta),\eta)$ are corresponded to the Hamilton
flow of $P$. Hence by (\ref{solveeikonal}) and a similar formula
for $\phi_P$, we have (\ref{phiphi}). This also shows that we can
solve (\ref{eikonalH}) in $(-t_0,t_0) \times \{x^2+\eta^2
<C\hbar^{1-\tau}\delta\}$.

To find $b_H$ we assume it is of the form (\ref{bH}) and we search
for functions $b_{H,j}$ such that
$e^{-i\phi_H/\hbar}(\frac{\hbar}{i}\d_t+\hat
H)(e^{i\phi_H/\hbar}b_H)\sim 0$.  After some straightforward
calculations and using the eikonal equation for $\phi_H$ we obtain
the so called transport equations (\ref{transportH}). Now let us
solve the transport equations inductively (see \cite{Ch}).

In \cite{Ch} it is shown that the solutions to the transport
equation (\ref{transportH}) are given by

\begin{equation} \label{solutiontotransport}
\begin{array}{lll}
  b_{H,0}(t,x,y,\eta,\hbar)=J^{-\half}(t,x,\eta,\hbar)
b_{H,0}(0,z(t,x,\eta;\hbar),\eta;\hbar), y, \eta,\hbar) \\ \\
  b_{H,j}(t,x,y,\eta,\hbar)=J^{-\half}(t,x,\eta,\hbar)\Big(
b_{H,j}(0,z(t,x,\eta;\hbar),\eta;\hbar), y, \eta,\hbar) \\
\qquad \qquad \qquad \qquad  -\half \int_0^t
J^{\half}(s,x,\eta,\hbar)\Delta
b_{H,j-1}(s,x(s,z(t,x,\eta;\hbar),\eta;\hbar), y, \eta,\hbar)ds \Big). \\
\end{array}
\end{equation}
where
$$J(t,x,\eta,\hbar)=\text{det}(\d_x z(t,x,\eta;\hbar))^{-1}.$$
Now, we notice by the assumption on $\psi_\hbar$, we have
supp$(b_{H,j}(0,x,y,\eta;\hbar))\subset \{x^2+\eta^2, x^2+\eta^2
<4c^{-1}\hbar^{1-\tau}\delta\}$. So by our previous discussion on
$z(t,x,\eta,\hbar)$, we can argue inductively that for all
$t\in(-t_0,t_0)$, supp$(b_{H,j})\subset \{x^2+\eta^2, y^2+\eta^2
<C_1\hbar^{1-\tau}\delta\}$ for some constant $C_1$. Since
$b_{H,j}|_{t=0}\in S_{\delta_0}^0$, we can also see inductively
from (\ref{solutiontotransport}) that $b_{H,j}\in S_{\delta_0}^0$.
Finally, Borel's theorem produces a compactly supported amplitude
$b_H \in S_{\delta_0}^0$ from the compactly supported functions
$b_{H,j}\in S_{\delta_0}^0$. This finishes the proof of items
$1-3$ in Lemma \ref{UH}.
\\

Now we give proofs for items $4,5$ in Lemma \ref{UH}.
\\

By choosing $C>C_1$, equation (\ref{phi_Hphi_P}) is clearly true
from (\ref{phiphi}). Next we prove that equation (\ref{b_Hb_P})
holds. Using (\ref{aPj}) and (\ref{aHj}), and because $P$ and $H$
agree in the ball $\{x^2+\eta^2<\frac{1}{4}\hbar^{1-2\ep}\}$, we
observe that the functions $a_{P,j}(x,\eta)$ and
$a_{H,j}(x,\xi,\hbar)$ agree in this ball. Therefore, because
supp$\psi_\hbar(x,\eta) \subset
\{x^2+\eta^2<4c^{-1}\hbar^{1-\tau}\delta\}$ and $\psi_\hbar=1$ in
$\{x^2+\eta^2<c\hbar^{1-\tau}\delta\}$, by (\ref{transportP}) and
(\ref{transportH})

$$b_{H,j}(0,x,y,\eta,\hbar)=b_{P,j}(0,x,y,\eta) \qquad \text{on} \quad
\{(x,y,\eta);\, x^2+\eta^2, \,
y^2+\eta^2<c\hbar^{1-\tau}\delta\}.$$ This proves (\ref{b_Hb_P})
only at $t=0$. But by applying (\ref{zPzH}) to
(\ref{solutiontotransport}) and a similar formula for $b_P$, we
get (\ref{b_Hb_P}). This finishes the proof of Lemma \ref{UH}.
\\

To finish the proof of Lemma \ref{criticallem}, we have to show
that for $t$ sufficiently small $Tr U_H(t)=Tr
U_P(t)+O(\hbar^\infty)$, or equivalently

$$\int\int
e^{i(\phi_H(t,x,\eta,\hbar)-x.\eta)/\hbar}
b_H(t,x,x,\eta,\hbar)dxd\eta=\int\int
e^{i(\phi_P(t,x,\eta)-x.\eta)/\hbar} b_P(t,x,x,\eta,\hbar)dxd\eta
+\, O(\hbar^\infty).$$ By (\ref{phi_Hphi_P}), the phase function
$\phi_H$ of the double integral on the left hand side equals
$\phi_P$ on the support of the amplitude $b_H$, so $\phi_H$ is
independent of $\hbar$ in this domain. Now, if $t\in(0,t_0)$ where
$t_0$ is smaller than the smallest non-zero period of the flows of
$P$ and $H$ respectively in the energy balls
$$\{(x,\eta)|\,\,H(x,\eta)\leq \delta \hbar^{1-\tau}\}\, \subset \,\{(x,\eta)|\,\,P(x,\eta)\leq \delta\},$$ then for every such
$t$, $(x,\eta)=(0,0)$ is the only critical point of the phase
functions $\phi_H(t,x,\eta,\hbar)-x.\eta$ and
$\phi_P(t,x,\eta)-x.\eta$ in these energy balls.

Obviously both integrals in the equation above are convergent
because their amplitudes are compactly supported. But the question
is whether or not we can apply the stationary phase lemma to these
integrals around their unique non-degenerate critical points. By
Lemma \ref{UH} the phase functions $\phi_H$ and $\phi_P$ are
independent of $\hbar$ on the support of their corresponding
amplitudes. Hence $\phi_H, \phi_P \in S_0^0$ on supp\,$b_H$ and
supp\,$b_P$ respectively. On the other hand
$b_H(t,x,x,\eta,\hbar)\in S^0_{\delta_0}$, $\delta_0<\half$; and
$b_P(t,x,x,\eta,\hbar)\in S^0_0$. These facts can be used to get
the required estimates for the remainder term in the stationary
phase lemma (for an estimate for the remainder term of the
stationary phase lemma, see for example Proposition $5.2$ of
\cite{DSj}).

Finally, by (\ref{phi_Hphi_P}) and (\ref{b_Hb_P}) it is obvious
that the integrals above must have the same stationary phase
expansions.

\section{Appendix B}

In this appendix we prove Lemma $\ref{reductionlemma}$. In fact we
prove that if $\Theta_\hbar$ is given by (\ref{Thetah}) then in
the sense of tempered distributions \begin{equation}
\label{reductionlemma2} Tr(\Theta_\hbar(\hat
H)e^{\frac{-it}{\hbar}\hat H})=Tr(e^{\frac{-it}{\hbar}\hat
H})+O(\hbar^\infty).\end{equation} Proof of Lemma
\ref{reductionlemma} follows similarly.
\\

We will use the min-max principle.
\\

\textbf{Min-max principle.} Let $H$ be a self-adjoint operator
that is bounded from below, i.e. $H\geq cI$, with purely discrete
spectrum $\{E_j\}_{j=0}^{\infty}$. Then

\begin{equation}\label{minmax} E_j= \sup_{\phi_1,...,\phi_{n-1}}
\inf_{\begin{array}{c}
  \psi \in D(H);\|\psi\|=1 \\
  \psi \in \text{span}(\phi_1,...,\phi_{n-1})^\perp\\
\end{array}}\; (\psi,H\psi).
\end{equation}
\\

As before we put $\hat H =-\half\hbar^2 \Delta + V_\hbar(x)=-\half
\hbar^2 \Delta +\half \sum_{k=1}^{n}\omega_k^2 x_k^2+W_\hbar(x)$,
and $\hat H_0=-\half \hbar^2 \Delta +\half
\sum_{k=1}^{n}\omega_k^2 x_k^2$. Then if we let
$C=\|W_\hbar(x)\|_{L^{\infty}\big( \mathbb R^n \times
(0,h_0)\big)}$, we have

$$(\psi,\hat H_0\psi)-C \leq (\psi,\hat H\psi)\leq (\psi,\hat H_0\psi)+C,$$
and therefore by applying the min-max principle to the operators
$\hat H$ and $\hat H_0$ we get

\begin{equation} \label{eigeninequality}E_j^{0}(\hbar)-C \leq E_j(\hbar)\leq E_j^{0}(\hbar)+C. \end{equation}
Notice we have explicit formulas for the eigenvalues
$E_j^{0}(\hbar)$ of $\hat H_0$. They are given by the lattice
points in the first quadrant of $\mathbb R^n$. More precisely

$$ \sigma(\hat H_0)=\Big\{E^0_{\vec \gamma}(\hbar)=\hbar\sum_{k=1}^n
\om_k(\gamma_k+\half);\;\; \gamma_k \in \mathbb Z^{\geq
0}\Big\}.$$
Since in the sense of tempered distributions
$$Tr(\Theta_\hbar(\hat H)e^{\frac{-it}{\hbar}\hat
H})=Tr(\chi_{[0,\delta \hbar^{1-\tau}]}(\hat
H)e^{\frac{-it}{\hbar}\hat H})+O(\hbar^\infty); \qquad (\text{see
for example \cite{KD}}),
$$to prove ($\ref{reductionlemma2}$), it is clearly enough to
show that for every $\phi$ in $S(\mathbb R)$

$$\sum_{\{j; \, E_j(\hbar)>\delta \hbar^{1-\tau}\}}\hat
\phi(\frac{E_j(\hbar)}{\hbar})=O(\hbar^{\infty}).$$ Since $\hat
\phi$ is in $S(\mathbb R)$, for every $p\geq 0$ there exists a
constant $C_p$ such that

$$ |\hat \phi (x)|\leq C_p |x|^{-p}.$$ Hence by $(\ref{eigeninequality})$

$$ \phi(\frac{E_j(\hbar)}{\hbar})\leq C_p
\Big|\frac{E_j(\hbar)}{\hbar}\Big|^{-p}\leq C_p
\Big|\frac{E_j^{0}(\hbar)-C}{\hbar}\Big|^{-p}.$$ Again using
$(\ref{eigeninequality})$ and because
$C=\|W_\hbar(x)\|_{L^{\infty}\big( \mathbb R^n \times
(0,h_0)\big)}<A\hbar^{\frac{3}{2}-3\ep}<\frac{\delta}{4}\hbar^{1-\tau}$
we get

$$ \phi(\frac{E_j(\hbar)}{\hbar})\leq  C_p
(\frac{\delta \hbar^{1-\tau}-C}{\delta \hbar^{1-\tau}-2C})^p
\Big|\frac{E_j^{0}(\hbar)}{\hbar}\Big|^{-p}<2C_p\Big|\frac{E_j^{0}(\hbar)}{\hbar}\Big|^{-p},
\qquad \quad \text{for}\;\; E_j(\hbar)>\delta \hbar^{1-\tau}.$$
Now let $m$ be an arbitrary positive integer. So in order to prove
the lemma it is enough to find a uniform bound for

$$A(\hbar):=\hbar^{-m}\sum_{\{\vec \gamma;\;\; \sum \om_k (\gamma_k +\half)>\frac{\delta \hbar^{1-\tau}-C}{\hbar}\}}|\sum_{k=1}^n
\om_k(\gamma_k+\half)|^{-p}.$$ By applying the
geometric-arithmetic mean value inequality we get

$$ \!\!\!\!\!\!\!\!\!\!\!\!\!\!\!\!\!\!\!\!\!\!\!\!\!\!\!\!\!\!\!\!\! A(\hbar) \leq n^{-p}\hbar^{-m}\sum_{\{\vec \gamma;\;\; \sum \om_k (\gamma_k +\half)>\frac{\delta \hbar^{1-\tau}-C}{\hbar}\}}|\prod_{k=1}^n
\om_k(\gamma_k+\half)|^{-p}$$

$$\quad \qquad \qquad \qquad
 \qquad \leq n^{-p} \sum_{k=1}^n\Big\{\Big(\hbar^{-m}\sum_{\{\gamma_k \in \mathbb Z^{\geq 0};\;\om_k (\gamma_k
+\half)>\frac{\delta
\hbar^{1-\tau}-C}{n\hbar}\}}|\om_k(\gamma_k+\half)|^{-p}\Big)
\prod _{k'\neq k}\Big(
\sum_{\gamma_{k'}}|\om_{k'}(\gamma_{k'}+\half)|^{-p}
\Big)\Big\}.$$ We claim for $p$ large enough there is a uniform
bound for the sum on the right hand side of the above inequality.
It is clear that if $p\geq 2$ then the series
$\sum_{\gamma_{k'}}|\om_{k'}(\gamma_{k'}+\half)|^{-p}$ is
convergent. Also if for some $\gamma_k$ we have $\om_k (\gamma_k
+\half)>\frac{\delta \hbar^{1-\tau}-C}{n\hbar}$, then because
$C=O(\hbar^{\frac{3}{2}-3\ep})$, for $\hbar$ small enough we have
$\big(\om_k (\gamma_k
+\half)\big)^{1/\tau}>(\frac{\delta}{2n})^{1/\tau}\frac{1}{\hbar}$.
Thus

$$\sum_{\{\gamma_k\in \mathbb Z^{\geq 0};\;\om_k (\gamma_k
+\half)>\frac{\delta
\hbar^{1-\tau}-C}{n\hbar}\}}\hbar^{-m}|\om_k(\gamma_k+\half)|^{-p}\leq
(\frac{2n}{\delta})^{m/\tau}
\sum_{\gamma_k}|\om_k(\gamma_k+\half)|^{\frac{m}{\tau}-p}.$$ So if
we choose $p>\max{\{\frac{m}{\tau},2\}}$, then the sum on the
right hand side is convergent and therefore we have a uniform
bound for the sum on the left hand side and hence for $A(\hbar)$.
This finishes the proof of (\ref{reductionlemma2}).

\end{document}